\documentclass[12pt]{amsart}
\usepackage{amssymb,amscd}
\usepackage{pstricks}

\newtheorem{theorem}{Theorem}[section]

\newtheorem{proposition}[theorem]{Proposition}

\newtheorem{definition}[theorem]{Definition}
\newtheorem{conj}[theorem]{Conjecture}

\theoremstyle{definition}
\newtheorem{remark}[theorem]{Remark}
\newtheorem{example}[theorem]{Example}
\numberwithin{equation}{section}

\newcommand{\Z}{\mathbb{Z}}
\newcommand{\R}{\mathbb{R}}
\newcommand{\C}{\mathbb{C}}
\newcommand{\CP}{\mathbb{CP}}
\newcommand{\RP}{\mathbb{RP}}

\renewcommand{\Re}{\mathrm{Re}\,}
\renewcommand{\Im}{\mathrm{Im}\,}
\newcommand{\hol}{\mathrm{hol}}

\title[Special Lagrangians, mirror symmetry and Calabi-Yau double covers]
{Special Lagrangian fibrations, mirror symmetry and Calabi-Yau double covers}
\author{Denis Auroux}
\thanks{Partially supported by NSF grants DMS-0600148 and DMS-0652630.}
\address{Department of Mathematics, M.I.T., Cambridge MA 02139, USA}
\email{auroux@math.mit.edu}

\begin{document}

\begin{abstract}
The first part of this paper is a review of the Strominger-Yau-Zaslow 
conjecture in various settings. In particular, we summarize how, given
a pair $(X,D)$ consisting of a K\"ahler manifold and an anticanonical
divisor, families of special Lagrangian tori in $X\setminus D$ and
weighted counts of holomorphic discs in $X$ can be used to build a
Landau-Ginzburg model mirror to $X$. In the second part
we turn to more speculative considerations about Calabi-Yau manifolds
with holomorphic involutions and their quotients. Namely,
given a hypersurface $H$ representing twice the anticanonical class in a
K\"ahler manifold $X$, we attempt to relate special
Lagrangian fibrations on $X\setminus H$ and on the (Calabi-Yau) double
cover of $X$ branched along $H$; unfortunately, the implications for
mirror symmetry are far from clear.
\end{abstract}

\maketitle

\vskip1.8cm

\section{Introduction}

The phenomenon of mirror symmetry was first evidenced for Calabi-Yau
manifolds, i.e.\ K\"ahler manifolds with holomorphically trivial canonical bundle. 
Subsequently it became apparent that mirror symmetry also holds in
a more general setting, if one enlarges the class of objects under
consideration (see e.g.\ \cite{HV}); namely, one should allow the mirror to be
a {\it Landau-Ginzburg model}, i.e.\ a pair consisting of a non-compact
K\"ahler manifold and a holomorphic function on it (called
{\it superpotential}).

Our motivation here is to understand how to {\it construct} the mirror
manifold, starting from examples where the answer is known and
extrapolating to less familiar situations; generally speaking,
the verification of the mirror symmetry conjectures for the
manifolds obtained by these constructions falls outside the scope of this
paper.

The geometric understanding of mirror symmetry in the Calabi-Yau case
relies on the Strominger-Yau-Zaslow (SYZ) conjecture \cite{SYZ}, which
roughly speaking postulates that mirror pairs of Calabi-Yau manifolds carry dual
fibrations by special Lagrangian
tori, and its subsequent refinements (see e.g.\ \cite{GS,KS}). This
program can be extended to the non Calabi-Yau case, as suggested by
Hori \cite{hori} and further investigated in \cite{Au}. In that case,
the input consists of a pair $(X,D)$ where $X$ is a compact K\"ahler
manifold and $D$ is a complex hypersurface representing the anticanonical
class. Observing that the complement of $D$ carries a holomorphic
$n$-form with poles along $D$, we can think of $X\setminus D$ as an
open Calabi-Yau manifold, to which one can apply the SYZ program.
Hence, one can attempt to construct the mirror of $X$
as a (complexified) moduli space of special Lagrangian tori in $X\setminus
D$, equipped with a Landau-Ginzburg superpotential defined by a weighted
count of holomorphic discs in $X$. However, exceptional discs and wall-crossing
phenomena require the incorporation of ``instanton corrections'' into the
geometry of the mirror (see \cite{Au}). 

One notable feature of the construction is that it provides a bridge between
mirror symmetry for the K\"ahler manifold $X$ and for the Calabi-Yau
hypersurface $D\subset X$. Namely, the fiber of the Landau-Ginzburg
superpotential is expected to be the SYZ mirror to $D$,
and the two pictures of homological mirror symmetry (for $X$ and for $D$)
should be related via restriction functors (see Section 7 of \cite{Au} for
a sketch). 
\medskip

In this paper, we would like to consider a slightly different
situation, which should provide another relation with mirror symmetry for
Calabi-Yau manifolds. The union of two copies of $X$ glued together along
$D$ can be thought of as a singular Calabi-Yau manifold, which can be smoothed to
a double cover of $X$ branched along a hypersurface $H$ representing twice
the anticanonical class and contained in a neighborhood of $D$. This
suggests that one might be able to think
of mirror symmetry for $X$ as a $\Z/2$-invariant version of mirror symmetry
for the Calabi-Yau manifold $Y$. Unfortunately, this proposal comes with
several caveats which make it difficult to implement.

Let $(X,\omega,J)$ be a compact K\"ahler manifold, and let $H$ be a
complex hypersurface in $X$ representing twice the anticanonical class.
Then the complement of $H$ carries a nonvanishing section $\Theta$ of
$K_X^{\otimes 2}$ with poles along $H$. We can think of $\Theta$
as the square of a holomorphic volume form defined up to sign.
In this context, we can look for special Lagrangian submanifolds
of $X\setminus H$, i.e.\ Lagrangian submanifolds on which the restriction
of $\Theta$ is real.
The philosophy of the SYZ conjecture suggests that, in favorable
cases, one might be able to construct a foliation of $X\setminus H$ in which
the generic leaves are special Lagrangian tori. Indeed, denote by $Y$ the
double cover of $X$ branched along $H$: then $Y$ is a Calabi-Yau manifold
with a holomorphic involution. If $Y$ carries a special Lagrangian
fibration that is invariant under the involution, then by quotienting we
could hope to obtain the desired foliation on $X\setminus H$; unfortunately
the situation is complicated by technicalities involving the symplectic
form.

\begin{conj}\label{conj:1} For a suitable choice of $H$,
$X\setminus H$ carries a special Lagrangian foliation whose lift to the
Calabi-Yau double cover $Y$ can be perturbed to a $\Z/2$-invariant
special Lagrangian torus fibration.
\end{conj}

If $-K_X$ is effective, we can consider
a situation where $H$ degenerates to a hypersurface $D$ representing
the anticanonical class in $X$, with multiplicity 2. As explained above, this corresponds to
the situation where $Y$ degenerates to the union of two copies of $X$ glued
together along $D$. One could hope that under such a degeneration
the foliation on $X\setminus H$ converges to a special Lagrangian torus
fibration on $X\setminus D$. Using the mirror construction described
in \cite{Au}, one can then try to relate a Landau-Ginzburg mirror
$(X^\vee,W)$ of $X$ to a Calabi-Yau mirror $Y^\vee$ of $Y$. The
simplest case should be when $K_{X|D}$ is holomorphically trivial (which
in particular requires $c_1(X)^2=0$).
Then $W:X^\vee\to \C$ is expected to have trivial monodromy around infinity
(see Remark \ref{rmk:relative}), so that
$\partial X^\vee\approx S^1\times D^\vee$ where $D^\vee$ is mirror to $D$.
It is then tempting to conjecture that, considering only the complex
structure of the mirror (and ignoring its symplectic geometry),
$Y^\vee$ can be obtained by gluing
together two copies of the mirror $X^\vee$ to $X$ along their boundary
$S^1\times D^\vee$. Unfortunately, as we will see in \S \ref{ss:mirror}
this is not compatible with instanton corrections.

The rest of this paper is organized as follows. In Section \ref{s:review}
we review the geometry of mirror symmetry from the perspective of the
SYZ conjecture, both in the Calabi-Yau case and in 
the more general case (relatively to an anticanonical divisor).
We then turn to more speculative considerations in
Section \ref{s:slaginvol}, where we discuss the geometry of Calabi-Yau
double covers, clarify the statement of Conjecture \ref{conj:1}, and
consider various examples. 

\subsection*{Acknowledgements}
I would like to thank Mohammed Abouzaid, Paul Seidel, Ludmil
Katzarkov, and Dima Orlov for many fruitful discussions.
I am also grateful to Ron Donagi,
whose interest in this topic prompted the writing of this paper. This
research was partially supported by NSF grants DMS-0600148 and DMS-0652630.

\section{The SYZ conjecture and mirror symmetry}\label{s:review}

\subsection{Motivation}\label{ss:motivation}
One of the most spectacular mathematical predictions of string theory
is the phenomenon of mirror symmetry, i.e.\ the existence of a broad
dictionary under which the symplectic geometry of a given manifold $X$
can be understood in terms of the complex geometry of a mirror manifold
$X^\vee$,
and vice-versa. This dictionary works at several levels, among which
perhaps the most exciting is Kontsevich's homological mirror conjecture,
which states that the derived Fukaya category of $X$ should be
equivalent to the derived category of coherent sheaves of its mirror
$X^\vee$
\cite{KoICM}; in the non Calabi-Yau case the categories under consideration
need to be modified appropriately \cite{KoENS} (see also \cite{Ab,HIV,KL,
SeVCM,SeBook}).

The main goal of the Strominger-Yau-Zaslow conjecture \cite{SYZ} is to
provide a geometric interpretation of mirror symmetry. Roughly speaking
it says that mirror manifolds carry dual fibrations by special Lagrangian
tori. In the Calabi-Yau case, one way to motivate the conjecture is to observe that, given any
point $p$ of the mirror $X^\vee$, mirror symmetry should put the
skyscraper sheaf $\mathcal{O}_p$ in correspondence with
some object $\mathcal{L}_p$ of the Fukaya category of
$X$. As a graded vector space 
$\mathrm{Ext}^*(\mathcal{O}_p,\mathcal{O}_p)$
is isomorphic to the cohomology of $T^n$; therefore the most likely candidate
for $\mathcal{L}_p$ is a (special) Lagrangian torus in $X$, equipped with
a rank 1 unitary local system (a flat $U(1)$ bundle). This suggests that
one should try to construct $X^\vee$ as a moduli space of pairs
$(L,\nabla)$ where $L$ is a special Lagrangian torus in $X$ and $\nabla$
is a flat unitary connection on the trivial line bundle over $L$.
Since for each torus $L$ the moduli space of flat 
connections can be thought of as a dual torus, we arrive at the familiar
picture.

When $X$ is not Calabi-Yau but the anticanonical class $-K_X$ is effective,
we can still equip the complement of a hypersurface $D\in |-K_X|$ with a
holomorphic volume form, and thus consider special Lagrangian
tori in $X\setminus D$. However, in this case, holomorphic discs in
$X$ with boundary in $L$ cause Floer homology to be {\it obstructed}\/
in the sense of Fukaya-Oh-Ohta-Ono \cite{FO3}:
to each object $\mathcal{L}=(L,\nabla)$ we can associate an obstruction
$\mathfrak{m}_0(\mathcal{L})$, given by a weighted count of holomorphic
discs in $(X,L)$, and the Floer differential on
$CF^*(\mathcal{L},\mathcal{L}')$ squares to $\mathfrak{m}_0(\mathcal{L}')-
\mathfrak{m}_0(\mathcal{L})$. Moreover, even when
the Floer homology groups $HF^*(\mathcal{L},\mathcal{L})$ can still be
defined, they are often zero, so that $\mathcal{L}$ is a trivial object
of the Fukaya category. On the mirror side, these features of the theory
can be replicated by the introduction of a {\it Landau-Ginzburg
superpotential}, i.e.\ a
holomorphic function $W:X^\vee\to\C$. Without getting into details, $W$
can be thought of as an obstruction term for the B-model on $X^\vee$,
playing the same role as $\mathfrak{m}_0$ for the A-model on $X$.
In particular, a point of $X^\vee$ defines a nontrivial object of the
category of B-branes $D^b_{sing}(X^\vee,W)$ only if it is a critical
point of $W$ \cite{KL,orlov}.

\subsection{Special Lagrangian fibrations and T-duality}\label{ss:t-dual}

Let $(X,\omega,J)$ be a smooth compact K\"ahler manifold of complex
dimension $n$. If $X$ is Calabi-Yau, i.e.\ the canonical bundle
$K_X$ is holomorphically trivial, then $X$ carries a globally defined
holomorphic volume form $\Omega\in \Omega^{n,0}(X)$:
this is the classical setting for mirror symmetry.
Otherwise, assume that $K_X^{-1}$ admits a nontrivial holomorphic section
$\sigma$, vanishing along a hypersurface $D$. Typically we will assume
that $D$ is smooth, or with normal crossing singularities.
Then $\Omega=\sigma^{-1}$ is a nonvanishing holomorphic $(n,0)$-form 
over $X\setminus D$, with poles along $D$.

The
restriction of $\Omega$ to a Lagrangian submanifold $L\subset X\setminus D$
does not vanish, and can be expressed in the form $\Omega_{|L}=\psi\,\mathrm{vol}_g$,
where $\psi\in C^\infty(L,\C^*)$ and $\mathrm{vol}_g$ is the volume form
induced on $L$ by the K\"ahler metric $g=\omega(\cdot,J\cdot)$.

\begin{definition}
A Lagrangian submanifold $L\subset X\setminus D$ is {\em special Lagrangian}
if the argument of $\psi$ is constant.
\end{definition}

The value of the constant depends only on the homology class $[L]\in
H_n(X\setminus D,\Z)$, and we will usually arrange for it to be a multiple
of $\pi/2$. For simplicity, in the rest of this paragraph we will assume
that $\Omega_{|L}$ is a real multiple of $\mathrm{vol}_g$.

The following classical result is due to McLean \cite{mclean} (at least
when $|\psi|\equiv 1$, which is the case in the Calabi-Yau setting;
see \S 9 of \cite{joycenotes} or Proposition 2.5 
of \cite{Au} for the case where $|\psi|\neq 1$):

\begin{proposition}[McLean]\label{prop:mclean}
Infinitesimal special Lagrangian deformations of $L$ are in one to one
correspondence with cohomology classes in $H^1(L,\R)$. Moreover, the
deformations are unobstructed.
\end{proposition}

More precisely, a section of the normal bundle $v\in C^\infty(NL)$
determines a 1-form $\alpha=-\iota_v\omega\in\Omega^1(L,\R)$ and an
$(n-1)$-form $\beta=\iota_v\mathrm{Im}\,\Omega\in\Omega^{n-1}(L,\R)$.
These satisfy $\beta=\psi\,*_g\alpha$, and the deformation is
special Lagrangian if and only if $\alpha$ and $\beta$ are both closed.
Thus special Lagrangian deformations correspond to ``$\psi$-harmonic''
1-forms $-\iota_v\omega\in \mathcal{H}^1_\psi(L)=\{\alpha\in \Omega^1(L,\R)|\,
d\alpha=0,\ d^*(\psi\alpha)=0\}$ (recall $\psi\in C^\infty(L,\R_+)$ is the
ratio between the volume elements determined by $\Omega$ and $g$).

In particular, special Lagrangian tori occur in $n$-dimensional families,
giving a local fibration structure provided that nontrivial $\psi$-harmonic
1-forms have no zeroes.

The base $B$ of a special Lagrangian torus fibration carries two natural
affine structures, which we call ``symplectic'' and ``complex''.
The first one, which encodes the symplectic geometry of $X$,
is given by locally identifying $B$ with a domain
in $H^1(L,\R)$ (where $L\approx T^n$).
At the level of tangent spaces, the cohomology class of 
$-\iota_v\omega$ provides an identification of $TB$ with
$H^1(L,\R)$; integrating, the local affine coordinates on $B$ are the symplectic
areas swept by loops forming a basis of $H_1(L)$. The other affine
structure encodes the complex geometry of $X$, and locally identifies
$B$ with a domain in $H^{n-1}(L,\R)$. Namely, one uses the cohomology
class of $\iota_v\mathrm{Im}\,\Omega$ to identify
$TB$ with $H^{n-1}(L,\R)$, and the affine coordinates are obtained by
integrating $\mathrm{Im}\,\Omega$ over the $n$-chains swept by cycles forming
a basis of $H_{n-1}(L)$.

In practice, $B$ can usually be compactified
to a larger space $\bar{B}$ (with non-empty boundary
in the non Calabi-Yau case), by also considering singular special
Lagrangian submanifolds that arise as limits of degenerating families
of special Lagrangian tori; however the affine structures
are only defined on the open subset $B\subset\bar{B}$.
\medskip

Ignoring singular fibers and instanton corrections, the first candidate
for the mirror of $X$ is therefore a moduli space $M$ of pairs $(L,\nabla)$,
where $L$ is a special Lagrangian torus in $X$ (or $X\setminus D$) 
and $\nabla$ is a flat
$U(1)$ connection on the trivial line bundle over $L$ (up to gauge).
The local geometry of $M$ is well-understood \cite{hitchin, leung, gross,
Au}, and
in particular we have the following result (see e.g.\ \S 2 of \cite{Au}):

\begin{proposition}\label{prop:mirrgeom}
$M$ carries a natural integrable complex structure $J^\vee$ arising from
the identification $$T_{(L,\nabla)}M=\{(v,\alpha)\in C^\infty(NL)\oplus
\Omega^1(L,\R)\,|\,-\iota_v\omega+i\alpha\in \mathcal{H}^1_\psi(L)\otimes
\C\},$$ a holomorphic $n$-form
$$\Omega^\vee((v_1,\alpha_1),\dots,(v_n,\alpha_n))=\int_L
(-\iota_{v_1}\omega+i\alpha_1)\wedge\dots\wedge
(-\iota_{v_n}\omega+i\alpha_n),$$
and a compatible K\"ahler form
$$\omega^\vee((v_1,\alpha_1),(v_2,\alpha_2))=\int_L \alpha_2\wedge
\iota_{v_1}\Im\Omega-\alpha_1\wedge \iota_{v_2}\Im\Omega$$
(this formula for $\omega^\vee$ assumes that
$\int_L \Re\Omega$ has been suitably normalized).
\end{proposition}

The moduli space of pairs $M$ can be viewed as a complexification of the moduli space
of special Lagrangian submanifolds; forgetting the connection gives a
projection map $f^\vee$ from $M$ to the real moduli space $B$. The fibers of
this projection are easily checked to be special Lagrangian tori in
$(M,\omega^\vee,\Omega^\vee)$. 

The special Lagrangian fibrations $f:X\to \bar{B}$ (or rather, its
restriction to the open subset $f^{-1}(B)$) and $f^\vee:M\to B$ can
be viewed as fiberwise dual to each other. In particular, it is easily
checked that the affine structure induced on $B$ by the symplectic
geometry of $f^\vee$ coincides with that induced by
the complex geometry of $f$, and vice-versa. Giving priority to the
symplectic affine structure, we will often implicitly equip $B$ with
the affine structure induced by the symplectic geometry of $X$, and
denote by $B^\vee$ the same manifold equipped with the other affine
structure (induced by the complex geometry of $X$, or the symplectic
geometry of $M$).

Thus, the philosophy of the Strominger-Yau-Zaslow
conjecture is that, in first approximation (ignoring instanton
corrections),  mirror symmetry amounts simply to exchanging
the two affine structures on $B$.
However, in general it is not at all
obvious how to extend the picture to the compactification $\bar{B}$.
The reader is referred to \cite{SYZ}, \cite{gross}, \cite{leung} for more details
in the Calabi-Yau case, and to \cite{hori} and \cite{Au} for the non
Calabi-Yau case.

\subsection{Mirror symmetry for Calabi-Yau manifolds}\label{ss:cyreview}

Constructing a special Lagrangian fibration on a Calabi-Yau manifold
is in general a challenging task, but there are a few situations where
it can be done explicitly, for instance in the case of flat tori,
or for hyperk\"ahler manifolds. We give two well-known
examples.

\begin{example}[Elliptic curves]\label{ex:elliptic}
Consider an elliptic curve $E=\C/(\Z\oplus \tau\Z)$, where
$\tau=i\gamma\in i\R_+$, equipped with the holomorphic volume form
$\Omega=dz$ and a K\"ahler form $\omega$ such
that $\int_E\omega=\lambda\in\R_+$. (The reason why we assume $\tau$
to be pure imaginary is that for simplicity we are suppressing any
discussion of $B$-fields). Then the family of circles parallel to the
real axis $\{\Im(z)=c\}$ defines a special Lagrangian fibration on $E$,
with base $B\simeq S^1$. One easily checks that the length of $B$ with
respect to the affine metric is equal to $\lambda$ for the symplectic
affine structure, and $\gamma$ for the complex affine structure.
The mirror elliptic curve $E^\vee$ is obtained by exchanging the two
affine structures on $B$; accordingly, it has modular parameter $\tau^\vee=
i\lambda$ and symplectic area $\int_{E^\vee}\omega^\vee=\gamma$.
(The reader is referred to \cite{PZ} for a verification of homological
mirror symmetry for the mirror pair $E,E^\vee$.)
\end{example}

\begin{example}[K3 surfaces]\label{ex:K3}
In the case of K3 surfaces, special Lagrangian fibrations can be built
using hyperk\"ahler geometry.
Let $(X,J)$ be an elliptically fibered K3 surface, for example obtained as
the double cover of $\CP^1\times\CP^1$ branched along a suitably chosen
algebraic curve of bidegree $(4,4)$: composing the covering map with
projection to the first $\CP^1$ factor, we obtain an elliptic fibration
$f:X\to\CP^1$ with 24 nodal singular fibers. 
Equip $X$ with a Calabi-Yau metric $g$, and denote the corresponding
K\"ahler form by $\omega_J$. Denote by
$\Omega_J$ a holomorphic $(2,0)$-form on $X$, suitably normalized, and let
$\omega_K=\mathrm{Re}(\Omega_J)$ and $\omega_I=\mathrm{Im}(\Omega_J)$:
then $(\omega_I,\omega_J,\omega_K)$ is a hyperk\"ahler triple for the
metric $g$.
Now switch to the complex
structure $I=g^{-1}\omega_I$ determined by the K\"ahler form $\omega_I$,
and with respect to which $\Omega_I=\omega_J+i\omega_K$ is a holomorphic
volume form. Since the fibers of
$f:X\to\CP^1$ are calibrated by $\omega_J$, the map $f$ is a special
Lagrangian fibration on $(X,\omega_I,\Omega_I)$.

The affine structures on the base of $f$ are only defined away from the
singularities of the fibration. Thus the geometry of $(X,\omega_I,\Omega_I)$
is characterized
by a pair of affine structures on the open subset $B\simeq S^2\setminus 
\{24\ \mathrm{points}\}$ of $\bar{B}\simeq S^2$. The monodromies
of the two affine structures around each singular point are the transpose
of each other, and each individual monodromy is conjugate to the standard
matrix {\tiny $\begin{pmatrix}1\!&\!1\\0\!&\!1 \end{pmatrix}$}. 

The mirror
of $(X,\omega_I,\Omega_I)$ is another K3 surface, carrying a special
Lagrangian fibration whose base differs from $B$ by an exchange of the
two affine structures. In fact, under certain assumptions (e.g., existence
of a section) and for a specific choice of $[\omega_J]$, the mirror may
be obtained simply by performing another hyperk\"ahler rotation to get
$(X,-\omega_K,\Omega_{-K}=\omega_J+i\omega_I)$; see e.g.\ \S 7 of 
\cite{huybrechts}. The reader is also referred to \S 7 of \cite{gross} for
more details on the SYZ picture for K3 surfaces.
\end{example}

In the above examples, one can avoid confronting heads-on the delicate
issues that arise when trying to reconstruct the mirror from the
affine geometry of $B$. In general, however, the compactification of
the mirror fibration over the singularities of the affine structure
and the incorporation of instanton corrections are two extremely challenging
aspects of this approach. The reader is referred to \cite{KS} and \cite{GS}
for two attempts at tackling this problem.

Another even more important issue is constructing a special
Lagrangian torus fibration on $X$ in the first place. When there is
no direct geometric construction as in the above examples, the most
promising approach seems to be Gross and Siebert's program
to understand mirror symmetry via toric degenerations \cite{GS03,GS}.
The main idea is to degenerate $X$ to a union $X_0$ of toric varieties
glued together along toric strata; toric geometry then provides a
special Lagrangian fibration on $X_0$, whose base is a polyhedral
complex formed by the union of the moment polytopes for the components
of $X_0$. Gross and Siebert then analyze carefully the behavior of this
special Lagrangian fibration upon deforming $X_0$ back to a smooth manifold,
showing how to insert singularities into the affine structure to compensate
for the nontriviality of the normal bundles to the singular strata along
which the smoothing takes place.
Moreover, they also show that, in the toric degeneration limit,
exchanging the affine structures on the base of the special Lagrangian
fibration can be understood as a combinatorial process called {\it discrete
Legendre transform} \cite{GS03}.

\begin{remark}\label{rmk:rays}
The affine geometry of $B$ is a remarkably powerful tool to understand
the symplectic and complex geometry of $X$ (and, by exchanging the
affine structures, of its mirror). Namely, away from the singularities,
the two affine structures on $B=B^\vee$ each determine an integral lattice 
in the tangent bundle $TB$; denoting these lattices by $\Lambda$ for
the symplectic affine structure and $\Lambda^\vee$ for the complex
affine structure, locally $X$ can be identified with either one of the
torus bundles $T^*B/\Lambda^*$ (with its standard symplectic form) and
$TB^\vee/\Lambda^\vee$ (with its standard complex structure).
Thus, locally, an integral affine submanifold of $B$ (i.e., a submanifold
described by linear equations with integer coefficients in local affine
coordinates with respect to the symplectic affine structure) determines
a Lagrangian submanifold of $X$ by the conormal construction. Similarly,
an integral affine submanifold with respect to the complex affine
structure $B^\vee$ locally determines a complex submanifold of $X$
(by considering its tangent bundle). More generally, {\it tropical
subvarieties} of $B$ or $B^\vee$ determine piecewise smooth Lagrangian
or complex subvarieties in $X$; whether these can be smoothed is a
difficult problem whose answer is known only in simple cases.

To give a concrete example, let us return to K3 surfaces
(Example \ref{ex:K3}) and the corresponding affine structures on 
$B\simeq S^2\setminus \{24\ \mathrm{points}\}$. Each singular fiber
of the special Lagrangian torus fibration $f:X\to\CP^1$ 
has a nodal singularity obtained
by collapsing a circle in the smooth fiber. The homology class of this
vanishing cycle determines a pair of rays in $B$ (straight half-lines emanating
from the singular point), with the property that the conormal bundles to
these rays compactify nicely to Lagrangian discs in $X$ (possibly after
a suitable translation within the fibers). Similarly, the
nodal singularity determines a pair of rays in $B^\vee$ (different from
the previous ones), whose tangent bundles (again after a suitable translation)
compactify to holomorphic
discs in $X$. When two singularities of the affine structure lie in a
position such that the corresponding rays in $B$ (resp.\ in $B^\vee$)
align with each other (and assuming the translations in the fibers also
match), the line segment joining them in $B$ (resp.\
$B^\vee$) determines a Lagrangian sphere (resp.\ a rational curve with
normal bundle $\mathcal{O}(-2)$)
in $X$. In the mirror $X^\vee$ the same alignment produces a rational
$-2$-curve (resp.\ a Lagrangian sphere). In fact, using the hyperk\"ahler
structure on $X$ and remembering that the elliptic fibration $f$ is 
$J$-holomorphic, these spheres correspond to (special) Lagrangian spheres in
$(X,\omega_J)$ which arise from the matching path construction and 
additionally are calibrated by $\omega_K$ (resp.\ $\omega_I$).
\end{remark}

\subsection{Mirror symmetry in the complement of an anticanonical
divisor}

We now consider special Lagrangian torus fibrations in the
complement $X\setminus D$ of an anticanonical divisor $D$ in a K\"ahler
manifold $X$. We start with a very easy example to make the following
discussion more concrete:

\begin{example} \label{ex:CP1}
Let $X=\CP^1$, equipped with any K\"ahler form invariant under the
standard $S^1$-action. Equip the complement
of the anticanonical divisor $D=\{0,\infty\}$, namely $\CP^1\setminus
\{0,\infty\}=\C^*$, with the standard holomorphic volume form
$\Omega=dz/z$. It is easy to check that the circles $|z|=r$
are special Lagrangian (with phase $\pi/2$). Thus we have a
special Lagrangian fibration $f:\CP^1\setminus D\to B$, whose base
$B$ is homeomorphic to an interval. As seen above, $B$ carries two
affine structures. With respect to the symplectic affine
structure, the special Lagrangian fibration is simply the moment map
for the $S^1$-action on $\CP^1$ (up to a factor of $2\pi$).
Thus $B$ is an open interval of
length equal to the symplectic area of $\CP^1$,
and can be compactified by adding the end points of the interval, which
correspond to the $S^1$ fixed points, i.e.\ the points of $D$. On the other
hand, with respect to the complex affine structure, $B$ is
an infinite line: from this point of view, the special
Lagrangian fibration is given by the map $z\mapsto \log |z|$.
\end{example}

We can start building a mirror to $X$ by considering the dual
special Lagrangian torus fibration $M$ as in \S \ref{ss:t-dual}.
$M$ is a non-compact K\"ahler manifold and, after taking
instanton corrections into account, it is in fact the mirror to
the open Calabi-Yau manifold $X\setminus D$. Thus, some information
is missing from this description. As explained at the end of \S
\ref{ss:motivation}, adding in the divisor $D$ very much affects
the special Lagrangian tori $X\setminus D$ from a Floer-theoretic
point of view, and the natural way to account for the resulting
obstructions is to make the mirror a Landau-Ginzburg model by
introducing a superpotential $W:M\to\C$. 

Recall that a point of $M$ is a pair $(L,\nabla)$, where $L\subset
X\setminus D$ is a special Lagrangian torus, and $\nabla$ is a flat
connection on the trivial line bundle over $L$. Given a homotopy
class $\beta\in \pi_2(X,L)$, we can consider the moduli space of
holomorphic discs in $X$ with boundary on $L$, representing the
class $\beta$. The virtual dimension (over $\R$) of this moduli space
is $n-3+\mu(\beta)$, where $\mu(\beta)\in\Z$ is the Maslov
index; in our case, the Maslov index is twice the algebraic
intersection number $\beta\cdot [D]$ (see e.g.\ Lemma 3.1 of
\cite{Au}). When $\mu(\beta)=2$, in favorable cases we can define
a (virtual) count $n_\beta(L)$ of holomorphic discs in the class
$\beta$ whose boundary passes through a generic point $p\in L$, and
define
\begin{equation}\label{eq:W}
W(L,\nabla)=\!\!\sum_{\substack{\beta\in \pi_2(X,L)\\ \mu(\beta)=2}}\!\!
n_\beta(L)\,z_\beta(L,\nabla),\ \ \mbox{where}\ \ z_\beta(L,\nabla)=
\exp(-{\textstyle \int_\beta\omega})\,\hol_{\nabla}(\partial\beta).
\end{equation}
Thus, $W$ is a weighted count of holomorphic discs of Maslov index 2
with boundary in $L$, with weights determined by the symplectic area
of the disc and the holonomy of the connection $\nabla$ along its
boundary. 

For example, in the case of $\CP^1$ (Example \ref{ex:CP1}), each special
Lagrangian fiber separates $\CP^1$ into two discs, one containing $0$
and the other one containing $\infty$. The classes $\beta_1$ and
$\beta_2$ represented by these discs satisfy $\beta_1+\beta_2=[\CP^1]$,
and hence the corresponding weights satisfy $z_{\beta_1}z_{\beta_2}=
\exp(-\int_{\CP^1}\omega)$. One can check that
$n_{\beta_1}=n_{\beta_2}=1$, so that using
$z=z_{\beta_1}$ as coordinate on $M$ we obtain the well-known
formula for the superpotential, $W=z+e^{-\Lambda}\,z^{-1}$, where $\Lambda$ is the symplectic
area of $\CP^1$.

While the example of $\CP^1$ is straightforward, several warnings are in
order. First, unless $X$ is Fano the sum
(\ref{eq:W}) is not known to converge. More importantly, if $L$ bounds
non-constant holomorphic discs of Maslov index 0 (i.e., discs contained
in $X\setminus D$), then the counts $n_\beta(L)$ depend on auxiliary data,
such as the point $p\in L$ through which the discs are required to pass,
or an auxiliary Morse function on $L$. 

An easy calculation shows that the
weights $z_\beta$ are local holomorphic functions on $M$ (with respect to
the complex structure defined in Proposition \ref{prop:mirrgeom}), and once
all ambiguities are lifted the disc counts $n_\beta(L)$ are 
locally constant, so that $W$ is locally a holomorphic function on $M$.
However, Maslov index 0 discs determine ``walls'' in $M$, across
which the counts $n_\beta(L)$ jump and hence the quantity (\ref{eq:W})
presents discontinuities. In terms of the affine geometry of the base of
the special Lagrangian fibration, an important mechanism for the generation
of walls comes from the rays in $B^\vee$ (the base with its
complex affine structure) that emanate from the vanishing cycles at the
singular fibers of the special Lagrangian fibration: indeed, by definition
any special Lagrangian fiber that lies on such a ray bounds a holomorphic
disc in $X\setminus D$ (see Remark \ref{rmk:rays}).
Intersections between these ``primary'' walls then generate further walls
(which can be visualized as rigid tropical configurations in $B^\vee$).

Fukaya-Oh-Ohta-Ono's results \cite{FO3} imply that
the formulas for $W$ in adjacent chambers of $M$ differ by a holomorphic
substitution of variables (see also Proposition 3.9 in \cite{Au}).
The guiding principle that governs {\it instanton
corrections} is that the various chambers of $M$ should be glued to each
other not in the naive manner suggested by the geometry of $B$, but rather
via the holomorphic gluing maps that arise in the wall-crossing formulas.
Thus, the instanton-corrected mirror is precisely the analytic space on
which the weighted count (\ref{eq:W}) of holomorphic discs in $(X,L)$,
and more
generally the ``open Gromov-Witten invariants'' of $(X,L)$ (yet to be
defined in the most general setting), become single-valued quantities.
The reader is referred to \cite{KS} and \cite{GS} for more details on
instanton corrections (in the Calabi-Yau case, but the general case is
similar).

One final issue is that, according to Hori and Vafa \cite{HV},
the mirror obtained by T-duality needs to be enlarged.
The holomorphic volume form $\Omega$ has poles along $D$, which causes
$B$ equipped with the complex affine structure to have infinite diameter
(after adding in the singular fibers, $B^\vee$ is complete). On the
other hand, the fact that $\omega$ extends smoothly across $D$ means
that, with respect to the symplectic affine structure, $B$ has finite
diameter, and compactifies to a singular affine manifold with boundary.
The consequence is that, after exchanging the affine structures, the
K\"ahler metric on the mirror is complete but its complex structure 
is ``incomplete'': for instance, in Example \ref{ex:CP1} the mirror of
$\CP^1$ is naturally a bounded annulus (of modulus equal to the symplectic
area of $\CP^1$), rather than the expected $\C^*$. Hori and Vafa's
suggestion (assuming that $X$ is Fano) is to symplectically ``enlarge''
$X\setminus D$ by considering a family of K\"ahler forms
$(\omega_k)_{k\to\infty}$ obtained by symplectic inflation along $D$,
with the property that $[\omega_k]=[\omega]+k\,c_1(X)$,
and simultaneously rescaling the
superpotential by a factor of $e^k$ (see also \S 4.2 of \cite{Au}).
However, this ``renormalization'' procedure is definitely not desirable
in the geometric setting considered in Section~\ref{s:slaginvol},
so we do not consider it further.

We end here our discussion of the various delicate points that come up in
the construction of the mirror and its superpotential, and simply
refer the reader to \cite{Au} for more details. Instead, we return to
examples.

\begin{example}[Toric varieties]\label{ex:toric}
Let $(X,\omega,J)$ be a toric variety of complex dimension $n$, and
consider the toric anticanonical divisor $D$ (i.e., the divisor of points
where the $T^n$-action is not free). Recall that $X\setminus D$ is
biholomorphic to $(\C^*)^n$, and equip it with the holomorphic
$(n,0)$-form $\Omega=d\log z_1\wedge\dots\wedge d\log z_n$, which has
poles along $D$.
Then the orbits of the standard $T^n$-action define a special Lagrangian
fibration on $X\setminus D\simeq (\C^*)^n$. With respect to the symplectic
affine structure, the base $B$ of this fibration is the moment polytope for
$(X,\omega)$, or rather its interior, and the special Lagrangian fibration
is simply given by the moment map. On the other hand, 
the complex affine structure on $B$ naturally identifies it with $\R^n$;
from this point of view the special Lagrangian fibration is the Log map
$(z_1,\dots,z_n)\mapsto (\log |z_1|,\dots,\log |z_n|)$.

Exchanging the two affine structures, the mirror of $X$ is naturally a
bounded domain in $(\C^*)^n$ (the subset of points whose image under the
Log map lies in the moment polytope of $X$), equipped with a complete
K\"ahler metric and a superpotential $W$ defined by a Laurent polynomial
consisting of one term for each component of $D$. Details can be found
in \cite{cho-oh} and \cite{FO3b} (see also \S 4 of \cite{Au} for a brief
overview, and \cite{Ab} for a partial verification of homological mirror
symmetry). 
\end{example}

\begin{example}[$\CP^2$]\label{ex:CP2}
Consider $\CP^2$ equipped with the Fubini-Study K\"ahler form $\omega_0$. Let
$D\subset \CP^2$ be a smooth elliptic curve defined by a homogeneous
polynomial of degree 3, and let $\Omega$ be a holomorphic
volume form on $\CP^2$ with poles along $D$.

\addtocounter{theorem}{-1}
\begin{conj} \label{conj:cp2} $\CP^2\setminus D$ carries a
special Lagrangian torus fibration over the disc with (generically) three nodal
singular fibers.
\end{conj}

\noindent
Tentatively, the construction of this special Lagrangian fibration 
proceeds as follows. Start with the toric setting, i.e.\ equip $\CP^2$
with a holomorphic volume form with poles along the toric
anticanonical divisor $D_0$ consisting of the three coordinate
lines ($\Omega_0=dx\wedge dy/xy$ in an affine chart). As mentioned above, the
orbits of the standard $T^2$-action define a special Lagrangian fibration on
$(\C^*)^2=\CP^2\setminus D_0$; with respect to the symplectic affine
structure, the base $B_0$ of this fibration is the moment polytope for $\CP^2$, 
i.e.\ a triangle. Deforming this situation to the case of a holomorphic
volume form $\Omega'$ with poles along a smooth cubic curve $D'$ obtained by
smoothing out the three nodal points of $D_0$ modifies the structure of
the special Lagrangian fibration near the three toric fixed points.
A local model for what happens near each of these points is described in
\S 5 of \cite{Au}. Namely, if we replace $\Omega_0$ by $\Omega_\varepsilon=
dx\wedge dy/(xy-\varepsilon)$, then the complement of the anticanonical divisor
$D_\varepsilon$ 
formed by the conic $xy=\varepsilon$ and the line at infinity carries a special
Lagrangian torus fibration with one nodal singular fiber: the fibers are
formed by intersecting the level sets of the moment map for the $S^1$-action
$e^{i\theta}\cdot(x,y)=(e^{i\theta}x,e^{-i\theta}y)$ with the level sets of
the function $|xy-\varepsilon|^2$, and the singularity is at the origin
\cite{Au}. If $\varepsilon$ is small then this family is close to the toric
family away from the origin. Therefore, general considerations about
deformations of families of special Lagrangians suggest that, if the smooth
elliptic curve $D'$ lies in a sufficiently small neighborhood of $D_0$, then
$(\CP^2\setminus D',\omega_0,\Omega')$ carries a special Lagrangian
fibration with three nodal singular fibers. From the point of view of the
affine geometry of the base $B'$ of this fibration, the smoothing of each node of 
$D_0$ amounts to replacing a corner of the triangle $B_0$ by a singular point
in the interior of $B'$ (so that $B'$ is a singular affine manifold with
boundary but without corners).

The special Lagrangian fibers over points close to the boundary of $B'$
lie in a tubular neighborhood of $D'$, and collapse to closed loops in
$D'$ as one approaches the boundary. Thus their first homology group is generated by
a meridian $m$ (the boundary of a small disc that
intersects $D'$ transversely once) and by a longitude $\ell$ (a curve
that runs parallel to a closed loop on $D'$). The monodromy of the special
Lagrangian fibration along
$\partial B'$ fixes $m$, but because the normal bundle to $D'$ has degree 9
it maps $\ell$ to $\ell+9m$. Thus, in a suitable basis the monodromy
along the boundary of $B'$ can be
expressed by the matrix {\tiny $\begin{pmatrix}1&9\\0&1\end{pmatrix}$}
(see equation (7.2) in \cite{Au}). 

The general case, where the cubic curve $D$ is not necessarily 
close to the singular toric configuration $D_0\subset \CP^2$,
should follow from a suitable result on deformations of
two-dimensional special Lagrangian torus fibrations with nodal singularities.
(To our knowledge such a result hasn't been proved yet; however it should
follow from an explicit analysis of the deformations of the nodal
singularities and the implicit function theorem applied to the smooth
part of the fibration. In our case one also needs to control the behavior
of the fibration near the boundary of $B$.)

When constructing the mirror, the singular fibers create walls, which
require instanton corrections. In the case of a cubic $D'$ obtained by
a small deformation of the toric configuration $D_0$, the local model
for a single smoothing suggests that the walls run
parallel to the boundary of the base $B'$. In fact, the special
Lagrangian fibers which lie sufficiently far from $D'$ are
Floer-theoretically equivalent to standard product tori. Thus, in the
``main'' chamber the superpotential is given by the same formula as
in the toric case, $W=x+y+e^{-\Lambda}/xy$ in suitable coordinates
(where $\Lambda=\int_{\CP^1}\omega$); in the other chambers it is
given by some analytic continuation of this expression (see \S 5 of \cite{Au}
for an explicit formula in the case of smoothing a single node of $D_0$).
In fact, ignoring completeness issues (e.g., looking only at $|W|\ll 1$),
the overall effect of deforming $D_0$ to a smooth cubic curve on the
complex geometry of the Landau-Ginzburg mirror is expected to be a
fiberwise compactification. Simultaneously, the symplectic area of the
fiber of the Landau-Ginzburg model, which is infinite in the toric case, 
is expected to become finite and equal to the imaginary part of the
modular parameter of the elliptic curve $D'$ (see also \cite{AKO}).
\end{example}

\begin{example}[Rational elliptic surface]\label{ex:E1}
Let $X$ be a rational elliptic surface obtained by blowing up $\CP^2$ at
the nine base points of a pencil of cubics, equipped with a K\"ahler form
$\hat{\omega}$. Let $\hat{D}\subset X$ be a smooth elliptic fiber 
(the proper transform of a cubic of the pencil), and let $\hat{\Omega}$
be a holomorphic $(2,0)$-form on $X$ with poles along $\hat{D}$. We
expect:

\addtocounter{theorem}{-1}
\begin{conj} \label{conj:ratell1} $X\setminus \hat{D}$ carries a
special Lagrangian torus fibration over the disc with (generically) 12 nodal singular
fibers. The monodromy of the affine structure around each singularity
is conjugate to {\tiny $\begin{pmatrix}1\!&\!1\\0\!&\!1
\end{pmatrix}$}, and the monodromy along $\partial \hat{B}$ is trivial.
\end{conj}

\noindent 
The construction starts with $(\CP^2, D, \omega_0, \Omega)$, where 
$D\subset\CP^2$ is an elliptic curve and $\Omega$ is a holomorphic
(2,0)-form with poles along $D$, as in Example \ref{ex:CP2} above.
By Conjecture \ref{conj:cp2}, we expect $\CP^2\setminus D$ to carry
a special Lagrangian torus fibration with three nodal singular fibers.
Now we blow up $\CP^2$ at nine points on the cubic $D$, to obtain the
rational elliptic surface $X$. 
Pulling back $\Omega$ under the blowup map yields a
holomorphic (2,0)-form $\hat\Omega$ on $X$, with poles along an elliptic
curve $\hat{D}\subset X$ (the proper transform of $D$).
On the other hand,
the K\"ahler form $\hat\omega$ on $X$ is not canonical, and depends
in particular on the choice of the symplectic areas of the exceptional divisors.
We claim that, provided these areas are sufficiently small, the blowup
should carry a special Lagrangian torus fibration with 12 nodal singular
fibers. 

The local model for each blowup operation is as follows \cite{AAK}.
Consider a neighborhood of the origin in $\C^2$ equipped with the standard
symplectic form, the holomorphic volume form $dx\wedge dy/y$ with
poles along $\C\times\{0\}$, and the family of special Lagrangian
cylinders $\{\mathrm{Re}(x)=t_1,\ \frac12 |y|^2=t_2\}\subset\C\times\C^*$.
Equip the blowup $\hat{\C}^2$ with a toric K\"ahler form $\hat{\omega}_0$
(invariant under the standard $T^2$-action) for which the area
of the exceptional divisor is $\epsilon>0$, and the holomorphic volume form
$\hat{\Omega}_0$ obtained by pulling back $dx\wedge dy/y$ under the blowup map
$\pi:\hat{\C}^2\to\C^2$.
The lift to $\hat{\C}^2$ of the $S^1$-action 
$e^{i\theta}\cdot(x,y)=(x,e^{i\theta}y)$ preserves $\hat\omega_0$
and $\hat\Omega_0$; its fixed point set consists of on one hand the proper 
transform $\hat{D}_0$ of $\C\times\{0\}$, and on the other hand the point where 
the proper transform of $\{0\}\times\C$ hits the exceptional divisor.
Denote by $\mu:\hat{\C}^2\to \R$ the moment map for this $S^1$-action,
normalized to equal $0$ on $\hat{D}_0$ and
$\epsilon$ at the isolated fixed point. Then it is easy to check that
the submanifolds
$\{\mathrm{Re}(\pi^*x)=t_1,\ \mu=t_2\}\subset \hat{\C}^2\setminus \hat{D}_0$
are special Lagrangian with respect to $\hat\omega_0$ and $\hat\Omega_0$ \cite{AAK}.
This family of special Lagrangians presents one nodal singular fiber --
the fiber which corresponds to $(t_1,t_2)=(0,\epsilon)$ and passes
through the isolated $S^1$-fixed point. 
Moreover, if $\epsilon$ is small then away from a neighborhood of the
exceptional divisor this family is close to the initial family of special
Lagrangians in $\C\times\C^*$. 

Even though the local model is only an asymptotic description
of the geometry of the special Lagrangian fibration on $\CP^2\setminus D$
near a point of $D$,  it should be possible to glue this local
construction into the fibration of Conjecture \ref{conj:cp2}, and thereby construct a special
Lagrangian fibration on the rational elliptic surface $X$ obtained by
blowing up $\CP^2$ at 9 points on the elliptic curve $D$. Each blow-up
operation inserts a nodal singular fiber into the fibration; thus
the base $\hat{B}$ of the special Lagrangian fibration on $X$ presents 12 singular
points. From the point of view of the symplectic affine structure, an easy
calculation on the local model shows that each new singular point lies
at a distance from the boundary of $\hat{B}$ equal to the symplectic area of the
exceptional curve of the corresponding blowup; in fact the
exceptional curve can be seen as a complex ray that runs
from the singular point to the boundary of $\hat{B}$. Moreover, the
monodromy of the fibration along the boundary of $\hat{B}$ is trivial,
reflecting the fact that the anticanonical divisor $\hat{D}\subset X$ 
has trivial normal bundle.

The general case, where the 
exceptional divisors of the blowups are not assumed to have small symplectic
areas, should again follow from a careful analysis of deformations of
two-dimensional special Lagrangian torus fibrations with nodal singularities
(with the same caveats as in the case of $\CP^2$).

\end{example}

\begin{remark}\label{rmk:relative} Assume $D$ is smooth. Then
the holomorphic $(n,0)$-form $\Omega$ on $X\setminus D$ induces a
holomorphic volume form $\Omega_D=\mathrm{Res}_D(\Omega)$ on $D$: the
{\it residue} of $\Omega$ along $D$. It is reasonable to expect that,
as is the case in the various examples considered above, in a neighborhood
of $D$ the special Lagrangian fibration on $(X\setminus D,\omega,\Omega)$
consists of tori which are $S^1$-bundles over special Lagrangian
submanifolds of $(D,\omega_{|D},\Omega_D)$. As a toy example,
consider $X=D\times\C$, 
$\omega=\omega_D+\frac{i}{2}dz\wedge d\bar{z}$, and 
$\Omega=\Omega_D\wedge dz/z$: then the product any special Lagrangian
submanifold of $D$ with a circle centered at the origin in $\C$ is
easily seen to be special Lagrangian. We conjecture that the qualitative
behavior is the same in the general case; see \S 7 of \cite{Au} for
more details. 

Assuming that this picture holds, the special Lagrangian
fibration $f:X\setminus D\to B$ can be extended over the boundary
of $B$ by a special Lagrangian fibration on $D$. In particular, the
boundary of $B$, with the induced affine structures, is the base $B_D$
of an SYZ fibration on $D$. More precisely: with respect to
the symplectic affine structure, the compactified base $\bar{B}$ is a
singular affine manifold with boundary (and corners if $D$ has normal
crossings),
and its boundary is $B_D$. With respect to the complex affine structure,
$B^\vee$ (after adding in the interior singular fibers) is a complete
singular affine manifold, isomorphic to $\R_+\times B_D^\vee$ outside
of a compact subset.

As already seen in Example \ref{ex:CP2}, near $\partial B$ the monodromy
of the affine structures on $B$ is determined
explicitly by the affine structures on $B_D$ and by the first Chern class
of the normal bundle
to $D$. Indeed, given a fiber of $f$ near the boundary of $B$,
i.e.\ an $S^1$-fibered special Lagrangian $L\subset X\setminus D$, the
action of the monodromy on $H_1(L)$ can be determined by working in
a basis consisting of a meridian loop linking $D$ and $n-1$
longitudes running parallel to $D$; from this one deduces the corresponding
actions on $H^1(L)$ (monodromy of $B$) and $H^{n-1}(L)$ (monodromy of
$B^\vee$).

Next, we look at the mirror, and observe that near its boundary $M$
consists of pairs $(L,\nabla)$ where $L$ is an $S^1$-fibered special
Lagrangian contained in a neighborhood of $D$. Denote by $\delta\in
\pi_2(X,L)$ the homotopy class of a small meridian disc intersecting
$D$ transversely once (with boundary the meridian loop), and let
$z_\delta(L,\nabla)$ be the corresponding weight as in equation (\ref{eq:W}).
Then $z_\delta$ is a holomorphic function on $M$ near its boundary.
(In fact, $z_\delta$ is the
dominant term in the expression of the superpotential $W$ near $\partial M$,
as the meridian discs have the smallest symplectic area among all
Maslov index 2 holomorphic discs.) By construction, the boundary of $M$
corresponds to the case where the area of the meridian disc reaches zero,
i.e.\ $\partial M=\{|z_\delta|=1\}$. 

Consider the complex
hypersurface $M_D=\{z_\delta=1\}$ ($\subset \partial M$).
Geometrically, $M_D$ corresponds to
limits of sequences of pairs $(L,\nabla)$ where $L$ collapses onto a
special Lagrangian torus $\Lambda\subset D$ and the connection $\nabla$ has trivial
holonomy along the collapsed $S^1$-factor in $L$, i.e.\ is pulled back from a
flat connection on the trivial bundle over $\Lambda$. Thus $M_D$ is none other
than the SYZ mirror to $D$. Moreover, the restriction of $z_\delta$ to
$\partial M$ induces a locally trivial fibration $z_\delta:\partial M\to
S^1$ with fiber $M_D$. The monodromy of this fibration can be realized
geometrically as follows. Start with a pair $(L,\nabla)$ where $L$ is
almost collapsed onto $\Lambda\subset D$ and $\nabla$ has trivial holonomy along
the meridian loop (so $z_\delta\in\R_+$): then we can change the holonomy of
$\nabla$ along the meridian loop by adding to it a multiple of $\sigma^{-1}
\nabla\sigma$, where $\sigma$ is the defining section of $D$
and $\nabla$ is a suitable connection on $K_X^{-1}$. From there it
follows easily that the monodromy of the fibration $z_\delta:\partial M\to S^1$ is a
symplectomorphism of $M_D$ which geometrically realizes 
(as a fiberwise translation in the special Lagrangian fibration
$M_D\to B_D$ dual to the SYZ fibration on $D$)
the mirror
to the autoequivalence $-\otimes K_{X\,|D}^{-1}$ of $D^b\mathrm{Coh}(D)$.

This rich geometric picture naturally leads to a formulation of mirror
symmetry for the pairs $(X,D)$ and $(M,M_D)$; see \S 7 of \cite{Au} for
details.

\end{remark}

\section{Special Lagrangian fibrations and double covers}\label{s:slaginvol}

\subsection{Special Lagrangians and Calabi-Yau double covers}

Let $(X,\omega,J)$ be a smooth compact K\"ahler manifold of complex
dimension $n$, and let $s$ be a nontrivial holomorphic section of
$K_X^{-2}$. Unless otherwise specified we assume that the hypersurface
$H=s^{-1}(0)$ is smooth.
$\Theta=s^{-1}$ is a nonvanishing section of $K_X^{\otimes 2}$
over $X\setminus H$, with poles along $H$, and locally $\Omega=\Theta^{1/2}$
is a nonvanishing holomorphic $n$-form, defined up to sign. The
restriction of $\Theta$ to a Lagrangian submanifold $L\subset X\setminus H$
does not vanish, and can be expressed in the form $\eta\,\mathrm{vol}_g^2$,
where $\eta\in C^\infty(L,\C^*)$. By analogy
with the situation considered previously, we make the following
definition:

\begin{definition}
A Lagrangian submanifold $L\subset X\setminus H$ is {\em special Lagrangian}
if the argument of $\eta$ is constant. (In fact $\Theta$ will usually be normalized
so that $\eta$ is real).
\end{definition}

It is easy to see that, if $L\subset X\setminus H$ is an orientable
special Lagrangian submanifold, then over $L$ the holomorphic quadratic
differential $\Theta$ admits a globally defined square root $\Omega$.
Therefore Proposition \ref{prop:mclean} still applies in this setting;
since $\Omega_{|L}=\eta^{1/2}\mathrm{vol}_g$, special Lagrangian
deformations are now given by $\eta^{1/2}$-harmonic 1-forms on $L$.

As before, the base $B$ of a special Lagrangian torus fibration carries
two natural affine structures, one arising from the symplectic geometry
of $X$ and the other one arising from its complex geometry.
\medskip

We now turn to the Calabi-Yau double cover of $X$ branched along $H$,
namely the unique double cover $\pi:Y\to X$ with the property that
$\tilde\Theta=\pi^*\Theta$ admits a globally defined
square root $\tilde{\Omega}\in \Omega^{n,0}(Y)$. More explicitly,
the obstruction for $\Theta$
to admit a globally defined square root is given by an element of
$H^1(X\setminus H,\Z/2)\simeq \mathrm{Hom}(\pi_1(X\setminus H),\Z/2)$,
and we consider the branched cover with this monodromy. 

The complex geometry of $Y$ is fairly straightforward, as the complex
structure $\tilde{J}$ and the holomorphic volume form $\tilde{\Omega}$
are simply lifted from those of $X$ via $\pi$. In particular, it is easy
to check that $\tilde{\Omega}$ is well-behaved along the ramification
divisor. (To give the simplest example, consider the map
$z\mapsto z^2$ from $\C$ to itself: the pullback of 
$\Theta=z^{-1}dz^{\otimes 2}$ is $\tilde{\Theta}=4\, dz^{\otimes 2}$,
which has a well-defined square root $\tilde{\Omega}=2\, dz$.)

On the other hand, constructing a K\"ahler form on $Y$ requires some
choices, because the pullback form $\pi^*\omega$ is degenerate along
the ramification locus $\tilde{H}=\pi^{-1}(H)$. One approach is to view $Y$ as
a complex hypersurface in the total space of the line bundle $K_X^{-1}$
over $X$, equipped with a suitable K\"ahler metric. More directly,
one can equip $Y$ with a K\"ahler form $\tilde\omega=\pi^*\omega+
\epsilon\,\lambda$, where $\epsilon>0$ is a sufficiently small constant and $\lambda$ is an
exact real $(1,1)$-form whose restriction to the complex line $\mathrm{Ker}(d\pi)$
is positive at every point of the ramification locus. Any two forms obtained
in this manner are symplectically isotopic; for example one can
take $\lambda=-i\partial\bar\partial\phi$ where $\phi:Y\to [0,1]$ is
supported in a neighborhood of $\tilde{H}$, equal to $1$ on $\tilde{H}$,
and strictly concave in the normal directions at every point of $\tilde{H}$.

Thus, given a compact special Lagrangian submanifold $L\subset X\setminus H$, 
the two lifts of $L$ are in general not special Lagrangian submanifolds
of $Y$, even though the restriction of $\tilde\Omega$ has constant phase, 
because they are not necessarily Lagrangian for $\tilde\omega$. In very
specific cases (for instance in dimension 1 or in product situations)
this is not an issue, but in general one needs to deform the
lift of $L$ to a nearby special Lagrangian submanifold $\tilde{L}\subset Y$,
whose existence is guaranteed by the unobstructedness of deformations
(Proposition \ref{prop:mclean}) as long as
$\tilde\omega$ is sufficiently close to $\pi^*\omega$.

When considering not just one submanifold but a whole special Lagrangian
fibration on $X\setminus H$, it is natural to ask whether the lifts
can be similarly deformed to a special Lagrangian fibration on $Y$.
Away from $H$ and from the singular fibers,
we can rely on an implicit function theorem for special Lagrangian
fibrations which again follows from unobstructedness. In spite of the
wealth of results that have been obtained on singularities of special
Lagrangians and their deformations (see e.g.\ \cite{joycev}), to our
knowledge there is no general result that would yield a special
Lagrangian fibration on $Y$ from one on $X\setminus H$. Nonetheless,
it seems reasonable to expect that such a result might hold at least in
low dimensions if the K\"ahler form on $Y$ is chosen suitably and the family
of special Lagrangians only presents generic singularities.

Thus, Conjecture \ref{conj:1} can be stated more precisely as follows:

\begin{conj}\label{conj:1plus}\ 
\begin{enumerate}
\item $X$ carries a special Lagrangian fibration (or rather, foliation)
$f:X\to \bar{B}$, where $\bar{B}$ is a singular affine manifold with
boundary (with two affine structures), such that the
generic fibers of $f$ are special Lagrangian tori in $X\setminus H$,
and the fibers of $f$ above $\partial\bar{B}$ are special Lagrangians
with boundary in $H$.\medskip

\item $Y$ carries a special Lagrangian torus fibration $\tilde{f}:Y\to
\tilde{B}$, where $\tilde{B}$ is a singular affine manifold without boundary
(with two affine structures), obtained by gluing together two copies of
$\bar{B}$ along their boundary.
\end{enumerate}
\end{conj}

\noindent
Note that the boundaries of the two copies of $\bar{B}$ are identified using
the identity map, whereas the normal direction is reflected; thus this
is an orientation-reversing gluing, and the resulting singular affine manifold
$\tilde{B}$ admits an orientation-reversing involution whose fixed point
locus is the ``seam'' of the gluing.

\subsection{Example: $\CP^1$ and elliptic curves}\label{ss:dim1}

As our first example, we consider $X=\CP^1$ equipped with any K\"ahler form
and a holomorphic quadratic differential $\Theta$ with poles at a
subset $H\subset\CP^1$. 

We first consider the special case
$\Theta=dz^2/(z^2-a^2)$, with simple poles at $\pm a$ and a double pole
at infinity. Setting $a=0$, we recover
the classical situation discussed in Example \ref{ex:CP1}, in which
the circles centered at the origin are special Lagrangian.
For arbitrary $a$,
it follows from classical geometry that every
ellipse with foci $\pm a$ is special Lagrangian with phase $\pi/2$ for
$\Omega=\Theta^{1/2}=dz/\sqrt{z^2-a^2}$. Thus we get a special Lagrangian
foliation of $\C\setminus \{\pm a\}$ by this family of ellipses, the sole
noncompact leaf being the real interval $(-a,a)$.  The general case is
less explicit but essentially amounts to modifying the special Lagrangian
family in the same manner not only near zero but also near infinity.

More precisely, equip $\CP^1$ with a generic holomorphic
quadratic differential $\Theta=z^2\,dz^2/(z-a)(z-b)(z-c)(z-d)$
with poles at $H=\{a,b,c,d\}$. Then, for a suitable
choice of phase, $\CP^1\setminus H$ admits a special Lagrangian foliation
in which all the leaves are closed loops with the exception of two
noncompact leaves, each connecting two of the points of $H$ (say $a$ and
$b$ on one hand, and $c$ and $d$ on the other hand).  For instance, if
$a<b<c<d$ are real, then we have such a foliation (with phase $\pi/2$)
in which the two noncompact leaves are the real line segments
$(a,b)$ and $(c,d)$. Indeed, after removing the two intervals $[a,b]$ and
$[c,d]$, the quadratic differential $\Theta$ admits a well-defined square root
$\Omega$, which is a closed 1-form and hence has the same period (easily
checked to be pure imaginary) on any homotopically nontrivial embedded curve.
The general case follows from the same~argument.

From a symplectic point of view, the base $B$ of this foliation is again an
interval of length equal to the symplectic area of $\CP^1$. However, unlike
the situation of Example \ref{ex:CP1}, the affine structure induced on $B$ by the
holomorphic volume form identifies it with a {\it finite} interval: if
we normalize $\Omega$ so that the integral of $\mathrm{Im}\,\Omega$ over
each special Lagrangian fiber is 1, then the length of
this interval is equal to $\int_b^c \mathrm{Re}\,\Omega$.

The double cover of $\CP^1$ branched at $H$ is an elliptic curve $Y$, and the
family of special Lagrangians in $\CP^1\setminus H$ lifts to a smooth
special Lagrangian fibration on $Y$. The base $\tilde{B}\simeq S^1$ of
this fibration, and its two affine structures, are obtained by
doubling $B$ along its boundary. For instance, the symplectic area of $Y$
(which is the length of $\tilde{B}$ with respect to the symplectic affine
structure, cf.\ Example \ref{ex:elliptic}) is twice that of $\CP^1$, whereas the integral of
$\mathrm{Re}\,\tilde{\Omega}$ over a section of the special Lagrangian
fibration (which is the length of $\tilde{B}$ with respect to the complex
affine structure) is twice $\int_b^c \mathrm{Re}\,\Omega$.

\begin{remark}\label{rmk:finitebase}
With respect to the complex affine structure, the base $B$ of the special
Lagrangian foliation on $(\CP^1\setminus H,\Omega)$
is a finite interval, whereas the base $B_0$ of the special Lagrangian fibration
on $(\CP^1\setminus \{0,\infty\},\Omega_0=dz/z)$ has infinite size. The reason
is that, as $a,b\to 0$ and $c,d\to\infty$, the elliptic curve $Y$
degenerates to a curve with two nodal singularities, and the base
$\tilde{B}$ of its special Lagrangian fibration degenerates
to a union of two infinite intervals. On the other hand, the symplectic
structure on $Y$, which determines the length of the base with respect to the
other affine structure, is unaffected by the degeneration.
\end{remark}

\subsection{Example: Elliptic surfaces}\label{ss:ratell}

We revisit Example \ref{ex:E1}, and again denote by $X$ a rational elliptic
surface obtained by blowing up $\CP^2$ at the 9 base points of a
pencil of cubics, equipped with a K\"ahler form $\hat\omega$. We previously
considered a holomorphic volume form $\hat\Omega$ on $X$ with poles
along an elliptic fiber $\hat{D}$.
Now we equip $X$ with a section $\Theta$ of
$K_X^{\otimes 2}$, with poles along the union $H=D_+\cup D_-$ of two
smooth fibers of the elliptic fibration; for simplicity we assume that
$D_\pm$ lie close to a same smooth fiber $\hat{D}$, so that away from a
neighborhood of $\hat{D}$ the quadratic volume element $\Theta$ is close to
the square ${\hat\Omega}^{\otimes 2}$ of the volume form considered
in Example \ref{ex:E1}.

\begin{conj} \label{conj:ratell2}
The special Lagrangian fibration
on $X\setminus \hat{D}$ constructed in Conjecture \ref{conj:ratell1} deforms to a special
Lagrangian family on $X\setminus H$.
The base $B$ of this family is homeomorphic to a closed disc, and over its
interior the fibers are special Lagrangian tori, with the exception of
12 nodal singular fibers. The fibers
above $\partial B$ are
special Lagrangian annuli with one boundary component on $D_+$ and the
other on $D_-$.
\end{conj}

We now explain the geometric intuition behind this conjecture by considering
a simplified local model in which everything is explicit.
The actual geometry of $X$ near $\hat{D}$ differs from this local model by higher order
terms; however the local model is expected to accurately describe all
the qualitative features of the special Lagrangian families of
Conjectures \ref{conj:ratell1} and \ref{conj:ratell2} in a small
neighborhood of $\hat{D}$.

In a small neighborhood of the fiber $\hat{D}$, the elliptic fibration
$X\to\CP^1$ is topologically trivial, and even though it is not
holomorphically trivial, in first approximation we can consider a
local model of the form $E\times U$, where $E$ is an elliptic curve
($E\simeq \hat{D}$) and $U$ is a neighborhood
of the origin in $\C$ (with coordinate $z$). In this simplified local model,
the holomorphic volume form $\hat\Omega$ can be written in the form
$dw\wedge dz/z$, where $dw$ is a holomorphic 1-form on $E$
(in fact, the {\it residue} of $\hat\Omega$ along $\hat{D}$), the symplectic
form $\hat\omega$ is a product form, and the
special Lagrangian family of Conjecture \ref{conj:ratell1} consists of
product tori, where the first factor is a special Lagrangian circle in
$(E,dw)$ and the second factor is a circle centered at the origin.

We now equip $E\times U$ with the quadratic volume element $\Theta=(dw\wedge dz)^{\otimes 2}
/(z^2-\epsilon^2)$, with poles along $H=E\times \{\pm \epsilon\}$.
Then the previous family of special Lagrangians
deforms to one where each submanifold is again a product: the first
factor is still a special Lagrangian circle in $(E,dw)$, and the second
factor is now an ellipse with foci at $\pm \epsilon$ (in the degenerate
case, the line segment $[-\epsilon,\epsilon]$).

The bases of these two special Lagrangian fibrations on $E\times U$,
equipped with their symplectic affine structures, are naturally
isomorphic, as each ellipse with foci at $\pm\epsilon$ can be used
interchangeably with the circle that encloses the same symplectic area
(in fact, the corresponding product Lagrangian tori in $E\times U$
are Hamiltonian isotopic to each other). In this sense, passing from
$X\setminus \hat{D}$ to $X\setminus H$ (i.e., from $\hat{B}$ to $B$)
is expected to be a trivial operation from the symplectic point of view.
However, the complex affine structures on $\hat{B}$ and $B$
are very different: from that perspective $\hat{B}$ is ``complete''
(its boundary lies ``at infinity'', since the affine structure
blows up near $\partial \hat{B}$ due to the singular behavior of
$\hat\Omega$ along $\hat{D}$), whereas $B$ has finite diameter.
This is most easily seen in terms of the local model near $\hat{D}$,
which allows us to reduce to the one-dimensional case (see
Remark \ref{rmk:finitebase}).
\medskip

Finally, we consider the double cover $Y$ of the rational elliptic
surface $X$ branched along $H$. It is
easy to see that $Y$ is an elliptically fibered K3 surface,
carrying a holomorphic involution under which the holomorphic volume
form $\tilde\Omega=(\pi^*\Theta)^{1/2}$ is anti-invariant.
By Conjecture \ref{conj:1} we expect that $Y$, equipped with a suitably
chosen K\"ahler form in the class $[\pi^*\hat\omega]$, carries a special
Lagrangian fibration with 24 nodal singular fibers, whose base
$\tilde{B}\simeq S^2$
is obtained by doubling $B$ along its boundary.
\medskip

In fact, it is well-known that such a fibration can be readily
obtained using hyperk\"ahler geometry as in Example \ref{ex:K3}. Indeed, consider an
elliptically fibered K3 surface with a real structure for which the
real part consists of two tori. For example, let
$Y'$ be the double cover of $\CP^1\times\CP^1$ branched along the zero
set of a generic real homogeneous polynomial of bidegree $(4,4)$ without
any real roots. Composing the covering map with projection to the first
$\CP^1$ factor, we obtain an elliptic fibration $f:Y'\to\CP^1$ with
24 singular fibers. Complex conjugation lifts to an
involution $\iota$ on $Y'$ which is antiholomorphic with respect to the
given complex structure $J$, and whose fixed point locus is the trivial
(disconnected) double cover of $\RP^1\times \RP^1$ (i.e., two tori).
The involution $\iota$ maps each fiber of $f$ to the fiber above
the complex conjugate point of $\CP^1$, and in particular it
interchanges pairs of complex conjugate singular fibers.

Equip $Y'$ with a Calabi-Yau metric, such that the K\"ahler form $\omega_J$ is
anti-invariant under $\iota$ (this is guaranteed by uniqueness of the
Calabi-Yau metric if one imposes $[\omega_J]$ to be the pullback of a
K\"ahler class on $\CP^1\times\CP^1$ and hence anti-invariant). Denote by
$\Omega_J$ a holomorphic $(2,0)$-form on $Y'$: then $\iota^*\Omega_J$ is a
scalar multiple of $\bar\Omega_J$, because $\dim H^{0,2}_J(Y)=1$. So after
normalization we can assume that $\iota^*\Omega_J=-\bar\Omega_J$, i.e.\
$\omega_K:=\mathrm{Re}(\Omega_J)$ is anti-invariant and
$\omega_I:=\mathrm{Im}(\Omega_J)$ is invariant. 

Now switch to the complex
structure $I$ determined by the K\"ahler form $\omega_I$. Then $\iota$
becomes a holomorphic involution, and the holomorphic volume form
$\Omega_I=\omega_J+i\omega_K$ is anti-invariant. Since the fibers of
$f:Y'\to\CP^1$ are calibrated by $\omega_J$, the map $f$ is a special
Lagrangian fibration on $(Y',\omega_I,\Omega_I)$, compatible with the
involution $\iota$.

It seems likely that this construction can be used as an alternative
approach to Conjecture \ref{conj:ratell2}, by considering the quotient
of this special Lagrangian fibration by the involution $\iota$.

\begin{remark}
The elliptic surface $X$ contains nine exceptional spheres, arising from
the nine blow-ups performed on $\CP^2$; these spheres 
intersect $H$ in two points, so their
preimages in the double cover $Y$ are rational curves with normal bundle
$\mathcal{O}(-2)$. These curves can be seen by looking at the
complex affine structures on the bases $B$ and $\tilde{B}$ of the
special Lagrangian fibrations on $X$ and $Y$, as discussed in Remark
\ref{rmk:rays}. Namely, the exceptional curves in $X$ correspond to
complex rays that run from singularities of the affine structure
of $B$ to its boundary (as in Example \ref{ex:E1}). Doubling $B$ along
its boundary to form $\tilde{B}$ creates alignments between pairs
of singular points lying symmetrically across from each other.
For at least 9 of the 12 pairs of points (those which correspond to the
blowups)
the corresponding complex rays match up to yield $-2$-curves in $Y$.
\end{remark}

\subsection{Example: $\CP^2$ and K3}\label{ss:cp2k3}

We now revisit Example \ref{ex:CP2}, and now equip $\CP^2$ with a
section $\Theta$ of $K^{\otimes 2}$ with poles along a smooth curve $H$ of
degree 6. We assume that $H$ lies in a small neighborhood of a cubic $D$,
i.e.\ it is defined by a homogeneous polynomial of the form
$p=\sigma^2-\epsilon q$, where $\sigma\in H^0(\mathcal{O}(3))$ is the
defining section of $D$ and $\epsilon$ is a small constant. Thus, away
from a neighborhood of $D$ the quadratic volume element $\Theta$ is close
to the square $\Omega^{\otimes 2}$ of the volume form considered in Example
\ref{ex:CP2}.

\begin{conj} \label{conj:cp22}
The special Lagrangian fibration
on $\CP^2\setminus D$ constructed in Conjecture \ref{conj:cp2} deforms to
a special Lagrangian family on $\CP^2\setminus H$.
The base $B$ of this family is homeomorphic to a closed disc, and over its
interior the fibers are special Lagrangian tori, with the exception of
three nodal singular fibers. The fibers above $\partial B$ are
special Lagrangian annuli with boundary on $H$, with the exception of
18 pinched annuli (with one arc connecting the two boundaries collapsed
to a point).
\end{conj}

While we do not have a complete picture to propose, the rough idea is as
follows. Looking at the defining section $p=\sigma^2-\epsilon q$ of $H$,
away from the zeroes of $q$ we can think of $H$ as two parallel copies
of $D$, and special Lagrangians are expected to behave as in the
previous example. Namely, near $D$ a special Lagrangian in $\CP^2\setminus D$
looks like the product of a special Lagrangian $\Lambda\,(\simeq S^1)$ in
$D$ with a small circle
in the normal direction, and the corresponding special Lagrangian in
$\CP^2\setminus H$ should be obtained by replacing the circle factor by
a family of ellipses whose foci lie on $H$. In the degenerate limit case,
the ellipses become line segments joining the two foci, forming an annulus;
when $\Lambda$ passes through a zero of $q$, the corresponding line segment
is collapsed to a point, giving a pinched annulus.

In fact, we are unable to provide an explicit local model for this behavior
on $X\setminus H$. However, Conjecture \ref{conj:cp22} can be corroborated
by calculations on a local model for the double cover $Y$ of $X$ branched
along $H$.

Near a point of $D$, we can consider local coordinates $(u,v)$ on a domain
in $\C^2$ such that $D$ is defined by the equation $u=0$, and $H$ is defined
by the equation $u^2-\epsilon q(v)=0$ for some holomorphic function $q$. The
corresponding section of $K_X^{\otimes 2}$ is given by
$\Theta=(u^2-\epsilon q(v))^{-1}\,(du\wedge dv)^{\otimes 2}$.
As $\epsilon\to 0$, this converges to the square of the holomorphic
volume form $u^{-1}\,du\wedge dv$, for which the cylinders $\{\Re v=a,
\ |u|^2=r\}$ are special Lagrangians (the circle factor corresponds to the
direction normal to $D$, while the other factor corresponds to a local
model for a special Lagrangian in $D$).

In this local model the double cover of $\C^2$ branched along $H$ is the
hypersurface $Y\subset \C^3$ defined by the equation $z^2=u^2-\epsilon q(v)$.
The pullback of $\Theta$ under the projection map $(z,u,v)\mapsto (u,v)$
admits the square root $$\tilde{\Omega}=z^{-1}\,du\wedge dv=u^{-1}\,dz\wedge
dv.$$
It is worth noting that $\tilde{\Omega}$ is the natural holomorphic volume form
induced on $Y$ by the standard volume form of $\C^3$: 
denoting by $f=z^2-u^2+\epsilon q(v)$ the defining function of $Y$, we
have $df\wedge \tilde{\Omega}=dz\wedge du\wedge dv$. We equip $Y$ with
the restriction of the standard K\"ahler form
$\omega_0=\frac{i}{2}\,dz\wedge d\bar{z}+\frac{i}{2}\,du\wedge d\bar{u}+
\frac{i}{2}\,dv\wedge d\bar{v}$, which differs from the pullback of the
standard K\"ahler form of $\C^2$ by the extra term $\frac{i}{2}\,dz\wedge
d\bar{z}=\frac{i}{2}\partial\bar\partial |u^2-\epsilon q(v)|^2$. We claim
that the (possibly singular) submanifolds 
$$\tilde{L}_{a,b}=\{(z,u,v)\in Y\ |\ \Re(v)=a,\ \Re(u\bar{z})=
b\}\qquad (a,b)\in\R^2$$ are special Lagrangian with respect to 
$\tilde{\Omega}$ and $\omega_0$. Indeed, the vector field $\xi(z,u,v)=
(iu, iz, 0)$ is tangent to the submanifolds $\tilde{L}_{a,b}$, and the
1-forms $\iota_\xi \Im\tilde{\Omega}=\Re dv$ and $\iota_\xi\omega_0=
-d\,\mathrm{Re}(u\bar{z})+\frac{i}{2}dv\wedge d\bar{v}$ both vanish on $
\tilde{L}_{a,b}$.
Moreover, $\tilde{L}_{a,b}$ is singular if and only if it passes through
a point $(0,0,v_0)$ with
$v_0$ a root of $q$.

The involution $(z,u,v)\mapsto (-z,u,v)$ maps $\tilde{L}_{a,b}$ to
$\tilde{L}_{a,-b}$. Thus, the special
Lagrangian fibration $(z,u,v)\mapsto (\Re v,\Re(u\bar{z}))$ descends
to a family of submanifolds in $\C^2$, parameterized by the quotient of
$\R^2$ by the reflection $(a,b)\mapsto (a,-b)$, i.e.\ the closed
upper half-plane. The image of $\tilde{L}_{a,b}$ under this projection is
$$L_{a,b}=\{(u,v)\in\C^2\,|\,\Re(v)=a,\ \Re(\bar{u}\sqrt{u^2-\epsilon q(v)})=\pm
b\},$$
and behaves exactly as described above: fixing a value of $v$ (i.e., a point
of $D$), the intersection of $L_{a,b}$ with $\C\times\{v\}$ is an ellipse
with foci the two square roots of $\epsilon q(v)$ (i.e.\ the two points
where $H$ intersects $\C\times\{v\}$). For $b=0$ the ellipse degenerates to
a line segment; when $v$ is a root of $q$ the ellipses become circles and
the line segment collapses to a point. However, a quick calculation shows
that $L_{a,b}$ is not Lagrangian with respect to the standard K\"ahler form
on $\C^2$.

Thus, it may well be easier to construct a special
Lagrangian fibration on the double cover of $\CP^2$ branched at $H$
(namely, a K3 surface) than on $\CP^2\setminus H$. In fact, as in the
previous example, the easiest way to construct such a fibration is
probably through hyperk\"ahler geometry, starting from an elliptically
fibered K3 surface with a real structure for which the real part is a
smooth connected surface of genus 10. Let $P$ be a real
homogeneous polynomial of bidegree $(4,4)$ whose zero set in $\RP^1\times
\RP^1$ consists of nine homotopically trivial circles $C_1,\dots,C_9$
bounding mutually disjoint discs $D_i$, and let $Y'$ be the
double cover of $\CP^1\times\CP^1$ branched along the zero set of $P$
(over $\C$). Then complex conjugation lifts to a $J$-antiholomorphic
involution $\iota$ of $Y'$, whose fixed point locus is a connected
surface of genus 10, namely the preimage of $\RP^1\times \RP^1\setminus
(D_1\cup\dots\cup D_9)$ (whereas the fixed point set of the
composition of $\iota$ with the nontrivial deck transformation consists
of 9 spheres, the preimages of $D_1,\dots,D_9$).
After performing a hyperk\"ahler rotation as in \S \ref{ss:ratell}, we
obtain a new complex structure $I$ on $Y'$ with respect to which $\iota$
is holomorphic and the elliptic fibration induced by projection to a
$\CP^1$ factor is special Lagrangian.

\begin{remark}
The curve $H\subset \CP^2$ bounds a number of Lagrangian discs, arising
as relative vanishing cycles for degenerations of $H$ to a nodal curve.
For instance, considering a degeneration of $H$ to two intersecting cubics
singles out 9 such discs. The preimages of these discs are Lagrangian
spheres in the double cover $Y$, and can be seen by looking at the
symplectic affine structure on the bases $B$ and $\tilde{B}$ of the special
Lagrangian fibrations on $\CP^2$ and $Y$. Namely, $\tilde{B}$ is obtained
by doubling $B$ along its boundary, and 18 of its singular points are
aligned along the ``seam'' of this gluing. The rays emanating from these
singular points run along the seam, and match with each other to give
rise to Lagrangian spheres.
\end{remark}

\begin{remark}
Consider a singular K3 surface $Y_0$ with 9
ordinary double point singularities, obtained as the double cover
of $\CP^2$ branched along the union $H_0$ of
two intersecting cubics. The singularities of $Y_0$ can be either smoothed,
which amounts to smoothing $H_0$ to a smooth sextic curve, or blown up,
which is equivalent to blowing up $\CP^2$ at the 9 intersection points
between the two components of $H_0$. These two procedures yield respectively
the K3 surface considered in the above discussion, and the K3 surface
considered in \S \ref{ss:ratell}.
$Y_0$ admits a special Lagrangian fibration whose base $\tilde{B}_0$
presents 9 singularities with monodromy conjugate to
{\tiny $\begin{pmatrix}1\!&\!2\\0\!&\!1 \end{pmatrix}$}; viewing
$\tilde{B}_0$ as two copies of a disc glued along the boundary, these
9 singularities all lie along the seam of the gluing. Smoothing $Y_0$
replaces each ordinary double point by a Lagrangian sphere, and 
resolves the corresponding singularity of $\tilde{B}_0$ into a pair
of singular points aligned along the seam. Blowing up $Y_0$
replaces each ordinary double point by an exceptional curve, and
resolves the corresponding singularity of $\tilde{B}_0$ into a pair
of singular points lying symmetrically across from each other on
either side of the seam.
\end{remark}

\subsection{Towards mirror symmetry for double covers} \label{ss:mirror}

Conjecture \ref{conj:1plus} suggests
that a mirror $Y^\vee$ of the Calabi-Yau double cover $Y$ of $X$ branched
along $H$ can be
obtained by gluing two copies of the mirror of $X\setminus H$ along their
boundary. From the point of view of affine geometry, we start with
a special Lagrangian fibration $f^\vee:M\to B$ (T-dual to the special
Lagrangian fibration on $X\setminus H$), and glue together two copies
of $M$ using an orientation-reversing diffeomorphism of $\partial M$
which induces a reflection in each fiber of $f^\vee$ above $\partial B$. 

Arguably the ``usual'' mirror of $X$ arises from considering the complement
of an anticanonical divisor $D$, rather than the hypersurface $H$.
Consider a degeneration of $H$ under which it
collapses onto $D$ (with multiplicity~2). At the level of double covers,
this amounts to degenerating $Y$ to the union of two copies of $X$ glued
together along $D$. By Moser's theorem, this deformation affects the complex
geometry of $Y$ but not its symplectic geometry.
Hence, the special Lagrangian fibrations on $X\setminus H$ and $X\setminus
D$ can reasonably be expected to have the same base $B$, as long as we only consider the
symplectic affine structure. (The complex affine structures are
very different: in the case of $X\setminus D$ the complex affine structure
blows up near the boundary of $B$, while in the case of $X\setminus H$
it doesn't. See e.g.\ Remark \ref{rmk:finitebase}.)
So, as long as we only consider the complex geometry of
the mirror and not its symplectic structure, it should be possible to
construct the mirror of $Y$ simply by gluing two copies of the mirror of
$X\setminus D$ (which is also the mirror of $X$ without its superpotential).

A complication arises when the normal bundle to $D$ is
not holomorphically trivial. In that case, the family of special Lagrangians
in $X\setminus H$ presents additional singularities at the boundary of $B$;
these singularities are not directly visible in the special Lagrangian
fibration on $X\setminus D$. An example of this phenomenon is presented
in \S \ref{ss:cp2k3} (compare Conjecture \ref{conj:cp22} with Conjecture
\ref{conj:cp2}). Thus, when building $\tilde{B}$ out of two copies of
the base $B$ of the special Lagrangian fibration on $X\setminus D$, we
need to introduce extra singularities into the affine structure
along the seam of the gluing. This is essentially the same phenomenon
as in Gross and Siebert's program (where singularities of
the affine structure also arise from the nontriviality of the normal bundles
to the codimension 1 toric strata along which the smoothing takes place).

For simplicity, let us just consider the case where $D$ has trivial normal
bundle. In that case, the discussion in Remark \ref{rmk:relative} implies
that the boundary of the (uncorrected) mirror $M$ of $X\setminus D$ is
the product of $S^1$ with a complex hypersurface $M_D\subset \partial M$
(the uncorrected SYZ mirror to $D$). In fact, we have a trivial
fibration $z_\delta:\partial M\approx M_D\times S^1\to S^1$, where 
$z_\delta$ is the weight
associated to the homotopy class of a meridian disc (collapsing to
a point as the special Lagrangian torus $L$ collapses onto a special
Lagrangian submanifold of $D$, whence $|z_\delta|=1$ on $\partial M$).
The orientation-reversing diffeomorphism $\varphi:\partial M\to\partial M$
used to glue
the two copies of $M$ together corresponds to a reversal of the coordinate
dual to the class of the meridian loop. More precisely, view a point
of $\partial M$ as a pair $(\Lambda,\nabla)$ where $\Lambda$ is a special
Lagrangian torus in $D$ and $\nabla$ is a flat unitary connection on the
trivial bundle over $\Lambda\times S^1$ (here we use the triviality of
the normal bundle to $D$ to view nearby special Lagrangians in 
$X\setminus D$ as products $\Lambda\times S^1$ rather than $S^1$-bundles
over $\Lambda$). Then the gluing diffeomorphism $\varphi$ is given by
$(\Lambda,\nabla)\mapsto (\Lambda,\bar\nabla)$, where $\bar\nabla$ is the
pullback of $\nabla$ by the diffeomorphism $(p,e^{i\theta})\mapsto 
(p,e^{-i\theta})$ of $\Lambda\times S^1$. Thus, under the identification of
$\partial M$ with $M_D\times S^1$, the diffeomorphism $\varphi$ is the
product of the identity map in $M_D$ and the complex conjugation map
$z_\delta\mapsto \bar{z}_\delta=z_\delta^{-1}$ from $S^1$ to itself.

At this point it would be tempting to conclude that, 
if $K_{X|D}$ is holomorphically trivial, then a mirror of $Y$ can be
obtained (at least as a complex manifold) by gluing
together two copies of the mirror of $X$ along their boundary
$S^1\times M_D$, to obtain a Calabi-Yau variety with a holomorphic
involution given near the ``seam'' of the gluing by $z_\delta\mapsto
z_\delta^{-1}$. Unfortunately, in the presence of instanton corrections
this seems to always fail; in particular, the fibers of $z_\delta:\partial
X^\vee\to S^1$ above two complex conjugate points are not necessarily 
biholomorphic. The following example in complex dimension 2 (inspired
by calculations in \cite{AAK}) illustrates
a fairly general phenomenon.

\begin{example}
We consider again the local model for
blow-ups mentioned in Example \ref{ex:E1}, modified so the special
Lagrangian fibers are tori rather than cylinders \cite{AAK}.
Start with $\C^*\times \C$ equipped
the holomorphic volume form $d\log x\wedge d\log y$ with poles along
$\C^*\times\{0\}$, and blow up the point $(1,0)$ to obtain 
a complex manifold $X$ equipped the holomorphic volume form
$\Omega=\pi^*(d\log x\wedge d\log y)$, with poles along the proper
transform $D$ of $\C^*\times\{0\}$. Observe that the
$S^1$-action $e^{i\theta}\cdot(x,y)=
(x,e^{i\theta}y)$ lifts to $X$, and consider an $S^1$-invariant
K\"ahler form $\omega$ for which the area of the exceptional divisor is
$\epsilon$. Denote by
$\mu:X\to\R$ the moment map for the $S^1$-action, normalized to equal
0 on $D$ and $\epsilon$ at the isolated fixed point.
The $S^1$-invariant tori $L_{t_1,t_2}=\{\log|\pi^*x|=t_1,\ \mu=t_2\}$
define a special Lagrangian fibration on $X\setminus D$,
with one nodal singularity at the isolated fixed point (for
$(t_1,t_2)=(0,\epsilon)$) \cite{AAK}.

The base $B$ of this special Lagrangian fibration is a half-plane, with 
a singular point at distance $\epsilon$ from the boundary (and nontrivial
monodromy around the singularity),
as pictured in Figure \ref{fig:1}; we place the
cut above the singular point in order to better visualize wall-crossing
phenomena near the boundary of $B$. The complex rays emanating from
the singular point (one of which corresponds precisely to the exceptional
divisor of the blowup) are responsible for wall-crossing jumps in holomorphic
disc counts, and split the mirror $M$ into two chambers, which are
essentially the preimages of the left and right halves of the figure.

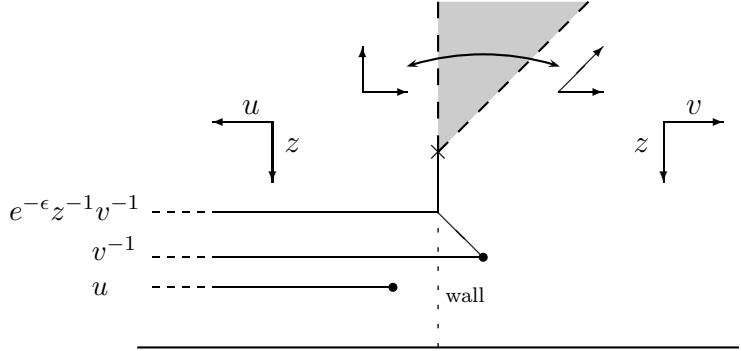
\begin{figure}[t]
\setlength{\unitlength}{2cm}
\begin{picture}(4.5,2.3)(-2.5,-0.3)
\psset{unit=\unitlength}
\newgray{ltgray}{0.8}
\pspolygon[fillstyle=solid,fillcolor=ltgray,linestyle=none](0,1)(0,2)(1,2)
\psline(-2,-0.3)(2,-0.3)
\put(0,1){\makebox(0,0)[cc]{\small $\times$}}
\psline[linestyle=dashed,dash=0.1 0.1](0,1)(0,2)
\psline[linestyle=dashed,dash=0.1 0.1](0,1)(1,2)
\psline[linestyle=dashed,dash=0.03 0.1,linewidth=0.5pt](0,-0.3)(0,1)
\put(0.05,0){\tiny wall}
\psarc{<->}(0.3,0){1.65}{72}{108}
\put(-0.5,1.4){\vector(1,0){0.3}}
\put(-0.5,1.4){\vector(0,1){0.3}}
\put(0.8,1.4){\vector(1,0){0.3}}
\put(0.8,1.4){\vector(1,1){0.3}}
\put(-1.1,1.2){\vector(-1,0){0.4}}
\put(-1.1,1.2){\vector(0,-1){0.4}}
\put(-1.02,1){$z$}
\put(-1.3,1.25){$u$}
\put(1.5,1.2){\vector(1,0){0.4}}
\put(1.5,1.2){\vector(0,-1){0.4}}
\put(1.3,1){$z$}
\put(1.65,1.25){$v$}
\pscircle[fillstyle=solid,fillcolor=black](0.3,0.3){0.03}
\pscircle[fillstyle=solid,fillcolor=black](-0.3,0.1){0.03}
\put(0.3,0.3){\line(-1,1){0.3}}
\put(0,0.6){\line(0,1){0.4}}
\put(0,0.6){\line(-1,0){1.5}}
\put(0.3,0.3){\line(-1,0){1.8}}
\put(-0.3,0.1){\line(-1,0){1.2}}
\multiput(-1.55,0.1)(-0.1,0){4}{\line(-1,0){0.05}}
\multiput(-1.55,0.6)(-0.1,0){4}{\line(-1,0){0.05}}
\multiput(-1.55,0.3)(-0.1,0){4}{\line(-1,0){0.05}}
\put(-2.3,0.28){\small $v^{-1}$}
\put(-2.85,0.55){\small $e^{-\epsilon}z^{-1}v^{-1}$}
\put(-2.3,0.05){\small $u$}
\end{picture}
\caption{A special Lagrangian fibration on the blowup of $\C^*\times\C$}
\label{fig:1}
\end{figure}

Denote by $z$ ($=z_\delta$) the holomorphic coordinate on $M$ which
corresponds to the holomorphic disc $\{\pi^*x=e^{t_1},\ \mu<t_2\}$ in
$(X,L_{t_1,t_2})$; it can be thought as a complexified and exponentiated
version of the downward-pointing affine coordinate pictured on Figure
\ref{fig:1}. In one of the two chambers of $M$, denote by $u$ the holomorphic
coordinate that similarly corresponds to the leftward-pointing
affine coordinate represented in the figure. For instance, if we
partially compactify $X$ to allow $\pi^*x$ to become zero (i.e., if we
had blown up $\C^2$ at $(1,0)$ rather than $\C^*\times\C$), then $u$
becomes (up to a scaling factor) the weight associated to a disc that runs
parallel to the $x$-axis. Similarly, denote by $v$ the holomorphic
coordinate in the other chamber of $M$ corresponding to a rightward-pointing
affine coordinate, normalized so that, if we ignore instanton corrections,
the gluing across the wall is given by $u=v^{-1}$.

Imagine that $L_{t_1,t_2}$ in the ``left'' chamber ($t_1<0$) bounds a
holomorphic disc
with associated weight $u$ (such a disc doesn't exist in $X$, but it exists
in a suitable partial compactification), and increase the value of $t_1$
past zero, keeping $t_2$ less than $\epsilon$: then this holomorphic disc
deforms appropriately (and its weight is now called $v^{-1}$), but it also
generates a new disc with weight $e^{-\epsilon}z^{-1}v^{-1}$, obtained
by attaching an exceptional disc (the part of the exceptional divisor where
$\mu>t_2$) as one crosses the wall. This phenomenon is pictured on Figure
\ref{fig:1} (where the various discs are abusively represented as
tropical curves, which actually should be drawn in the complex affine structure).
Thus the instanton-corrected gluing is
given by $u=v^{-1}+e^{-\epsilon}z^{-1}v^{-1}$, i.e.,
\begin{equation}\label{eq:corr} uv=1+e^{-\epsilon}z^{-1}.\end{equation}
Actually the portion of the wall where $t_2>\epsilon$ also gives rise to
the same instanton-corrected gluing, so that the corrected mirror is
globally given by (\ref{eq:corr}); see \cite{AAK}.

Now replace $D$ by the union $H=D_+\cup D_-$ of two disjoint complex curves, e.g.\
the proper transforms of two complex lines intersecting transversely
at the blown up point $(1,0)$, and consider the double cover $Y$ of $X$
branched along $H$. (We leave the details unspecified, as the
construction should arguably be carried out in a global setting such
as that of Conjecture \ref{conj:ratell2} rather than in the local setting.)

Conjecture \ref{conj:1} suggests that $Y$ should carry a special Lagrangian
fibration whose base (considering only the symplectic affine structure)
is obtained by doubling $B$ along its boundary. Pictorially, this
corresponds to flipping Figure \ref{fig:1} about the horizontal axis and
gluing the two pictures together. On the mirror, before instanton
corrections this amounts to reflecting the $z$ variable via $z\mapsto
z^{-1}$, and gluing $M$ and its reflected copy along their common boundary
$|z|=1$. However, the gluing via $z\mapsto z^{-1}$ is not compatible with
the instanton corrections discussed above; this is because when we
cross the wall there are now two different exceptional discs to consider.
Namely, $Y$ contains a $-2$-curve $C$ (the preimage of the
exceptional curve in $X$), corresponding to the alignment between the
walls that come out of the two singular fibers on either side of the seam.
Special Lagrangian fibers which lie on the wall intersect $C$ in a circle
and split it into two Maslov index 0 discs, which both contribute to
instanton corrections. A careful calculation shows that the
instanton-corrected gluing is now 
\begin{equation}\label{eq:corr2}
uv=(1+e^{-\epsilon}z^{-1})(1+e^{-\epsilon}z).
\end{equation}
Thus the instanton-corrected mirror to $Y$ does carry a holomorphic
involution defined by $z\mapsto z^{-1}$, but restricting to the subset $|z|<1$
does not yield the instanton-corrected mirror to $X$.
\end{example}


\begin{thebibliography}{99}
\bibitem{Ab}
M. Abouzaid,
{\sl Morse homology, tropical geometry, and homological mirror symmetry
for toric varieties}, math.SG/0610004.
\bibitem{AAK}
M. Abouzaid, D. Auroux, L. Katzarkov, in preparation.
\bibitem{Au}
D. Auroux,
{\sl Mirror symmetry and T-duality in the complement of an anticanonical
divisor},
J.~G\"okova Geom.\ Topol.\ {\bf 1} (2007), 51--91, arXiv:0706.3207.
\bibitem{AKO}
D. Auroux, L. Katzarkov, D. Orlov,
{\sl Mirror symmetry for Del Pezzo surfaces: Vanishing cycles and coherent
sheaves},
Inventiones Math. {\bf 166} (2006), 537--582, math.AG/0506166.
\bibitem{cho-oh}
C.-H. Cho, Y.-G. Oh,
{\sl Floer cohomology and disc instantons of Lagrangian torus fibers in
Fano toric manifolds}, Asian J. Math. {\bf 10} (2006), 773--814,
math.SG/0308225.
\bibitem{FO3}
K. Fukaya, Y.-G. Oh, H. Ohta, K. Ono,
{\sl Lagrangian intersection Floer theory: Anomaly and obstruction},
preprint, second expanded version, 2006.
\bibitem{FO3b}
K. Fukaya, Y.-G. Oh, H. Ohta, K. Ono,
{\sl Lagrangian Floer theory on compact toric manifolds I}, 
arXiv:0802.1703.
\bibitem{gross}
M. Gross,
{\sl Special Lagrangian Fibrations II: Geometry},
Winter School on Mirror Symmetry, Vector Bundles and Lagrangian Submanifolds
(Cambridge, MA, 1999), AMS/IP Stud.\ Adv.\ Math. {\bf 23}, Amer.\ Math.\
Soc., Providence, 2001, pp.\ 95--150,
math.AG/9809072.
\bibitem{GS03}
M. Gross, B. Siebert,
{\sl Affine manifolds, log structures, and mirror symmetry},
Turkish J.\ Math.\ {\bf 27} (2003), 33--60,
math.AG/0211094.
\bibitem{GS}
M. Gross, B. Siebert,
{\sl From real affine geometry to complex geometry},
math.AG/0703822.
\bibitem{hitchin}
N. Hitchin,
{\sl The moduli space of special Lagrangian submanifolds},
Ann.\ Scuola Norm.\ Sup.\ Pisa Cl.\ Sci.\ {\bf 25} (1997), 503--515,
dg-ga/9711002.
\bibitem{hori}
K. Hori,
{\sl Mirror symmetry and quantum geometry}, Proc.\ ICM (Beijing, 2002),
Higher Ed.\ Press, Beijing, 2002, vol.\ III, 431--443,
hep-th/0207068.
\bibitem{HIV}
K.\ Hori, A.\ Iqbal, C.\ Vafa, {\sl D-branes and mirror symmetry},
hep-th/0005247.
\bibitem{HV}
K. Hori, C. Vafa, 
{\sl Mirror symmetry}, hep-th/0002222.
\bibitem{huybrechts}
D. Huybrechts, {\sl Moduli spaces of hyperk\"ahler manifolds and mirror
symmetry}, Intersection theory and moduli, ICTP Lect.\ Notes {\bf 19}, ICTP,
Trieste, 2004, 185--247, math.AG/0210219.
\bibitem{joycenotes}
D. Joyce, {\sl Lectures on Calabi-Yau and special Lagrangian geometry},
math.DG/0108088.
\bibitem{joycev}
D. Joyce, {\sl Special Lagrangian submanifolds with isolated conical
singularities. V. Survey and applications}, J.\ Differential Geom.\ 
{\bf 63} (2003), 279--347, math.DG/0303272.
\bibitem{KL}
A.\ Kapustin, Y.\ Li, {\sl D-branes in Landau-Ginzburg models and
algebraic geometry}, J.\ High Energy Phys.\ {\bf 0312} (2003), 005,
hep-th/0210296.
\bibitem{KoICM}
M.\ Kontsevich, {\sl Homological algebra of mirror symmetry},
Proc.\ International Congress of Mathematicians (Z\"urich, 1994),
Birkh\"auser, Basel, 1995, pp.\ 120--139.
\bibitem{KoENS}
M.\ Kontsevich, {\sl Lectures at ENS, Paris, Spring 1998}, notes taken by
J.\ Bellaiche, J.-F.\ Dat, I.\ Marin, G.\ Racinet and H.\ Randriambololona,
unpublished.
\bibitem{KS}
M.\ Kontsevich, Y.\ Soibelman,
{\sl Affine structures and non-Archimedean analytic spaces},
The unity of mathematics, Progr. Math. {\bf 244}, Birkh\"auser Boston, 
2006, pp. 321--385, math.AG/ 0406564.
\bibitem{leung}
N.\,C. Leung, {\sl Mirror symmetry without corrections},
math.DG/0009235.
\bibitem{mclean}
R.\,C. McLean,
{\sl Deformations of calibrated submanifolds},
Comm. Anal. Geom. {\bf 6} (1998), 705--747.
\bibitem{orlov}
D.~Orlov, {\sl Triangulated categories of singularities and D-branes in
Landau-Ginzburg models}, math.AG/0302304.
\bibitem{PZ}
A. Polishchuk, E. Zaslow,
{\sl Categorical mirror symmetry: the elliptic curve},
Adv.\ Theor.\ Math.\ Phys.\ {\bf 2} (1998), 443--470.
\bibitem{SeVCM}
P. Seidel, {\sl Vanishing cycles and mutation},
Proc. 3rd European Congress of Mathematics (Barcelona, 2000), Vol. II,
Progr.\ Math.\ {\bf 202}, Birkh\"auser, Basel, 2001, pp.\ 65--85,
math.SG/0007115.
\bibitem{SeBook}
P. Seidel,
{\sl Fukaya categories and Picard-Lefschetz theory},
European Math.\ Soc., to appear.
\bibitem{SYZ}
A. Strominger, S.-T. Yau, E. Zaslow,
{\sl Mirror symmetry is T-duality},
Nucl. Phys. B {\bf 479} (1996), 243--259, hep-th/9606040.
\end{thebibliography}
\end{document}